\documentclass[11pt,letterpaper]{amsart}
\usepackage[utf8]{inputenc}
\usepackage{verbatim}
\usepackage{amsmath,amssymb,amsfonts,amsthm}
\usepackage{bigints}
\usepackage{hyperref}
\usepackage{array}

\usepackage{pdflscape}
\usepackage{multirow}
\usepackage{hyperref}

\usepackage{graphicx}
\usepackage{bbm}
\usepackage{mathtools}
\mathtoolsset{showonlyrefs}
\usepackage{amssymb}
\usepackage{verbatim}
\usepackage{amsmath,amssymb}
\usepackage{amsmath}

\newcommand{\norm}[1]{ \left\lVert#1\right\rVert}

\usepackage{bigints}
\usepackage{enumitem}

\usepackage{amsfonts,array,amsmath,xcolor}

\usepackage{amsthm}

\numberwithin{equation}{section}

\usepackage{graphicx}

\newtheorem{definition}{Definition}[section]
\newtheorem{proposition}[definition]{Proposition}
\newtheorem{theorem}[definition]{Theorem}

\newtheorem{lemma}[definition]{Lemma}
\newtheorem{corollary}[definition]{Corollary}
\newtheorem{remark}[definition]{Remark}

\begin{document}
\thanks{
This work is supported by
 PID2023-150984NB-I00 funded by MICIU/AEI/10.13039/501100011033/ FEDER, EU, 
the CERCA Programme of the Generalitat de Catalunya and the Severo Ochoa, and Mar\'ia de Maeztu
Program for Centers and Units of Excellence in R\&d (CEX2020-001084-M).
M. Saucedo is supported by  the Spanish Ministry of Universities through the FPU contract FPU21/04230.
 S. Tikhonov is supported
by 2021 SGR 00087. }

\author{Miquel Saucedo}
\address{M.  Saucedo,  Centre de Recerca Matemàtica\\
Campus de Bellaterra, Edifici C
08193 Bellaterra (Barcelona), Spain}
\email{miquelsaucedo98@gmail.com }

\author{Sergey Tikhonov}
\address{S. Tikhonov, 
ICREA, Pg. Lluís Companys 23, 08010 Barcelona, Spain;
Centre de Recerca Matem\`{a}tica\\
Campus de Bellaterra, Edifici C
08193 Bellaterra (Barcelona), Spain,
 and Universitat Autònoma de Barcelona, Spain.}
\email{stikhonov@crm.cat}

\subjclass[2010]{Primary  
42B10, 42B35; Secondary 46E30.}
\keywords{Fourier transform,
Uncertainty principles, Weighted Lebesgue spaces}

\title[Heisenberg's inequality in $L^p$] 
{Heisenberg's inequality in $L^p$}

\begin{abstract}
   In this paper, we obtain 
   non-symmetric and 
   symmetric
versions of the classical Heisenberg--Pauli--Weyl uncertainty principle
in Lebesgue spaces with power weights.

\end{abstract}
\maketitle
\section{Introduction}
Given an integrable function $f: \mathbb{R}^d\to \mathbb{C}$,
we define its Fourier transform by
$$ \widehat{f}(\xi)=\int_{\mathbb{R}^d} {f}(x) e^{-2\pi i \langle x, \xi \rangle} d x.
$$
We begin with the Heisenberg uncertainty principle (HUP), discovered almost one hundred years ago, which is nowadays one of the fundamental results in mathematics and  quantum
mechanics, 
 with numerous applications in both fields; see, for example, \cite{folland}. Heisenberg's estimate  states that if $f\in L_2(\mathbb{R}^d)$,  
then
\begin{equation}\label{heisenberg}
 \|f\|_2^4
\lesssim\int_{\mathbb{R}^d}
|x |^2
|f(x)|^2 dx
\int_{\mathbb{R}^d}
|\xi |^2 |\widehat{f}(\xi)|^2 d\xi,
\end{equation} or, equivalently,
$$\|f\|_2^2
\lesssim\int_{\mathbb{R}^d}
|x |^2
|f(x)|^2 dx +
\int_{\mathbb{R}^d}
|\xi |^2 |\widehat{f}(\xi)|^2 d\xi.$$

In this paper, we establish two extensions of 
Heisenberg's result in the setting of weighted Lebesgue spaces: a \textit{non-symmetric} and a \textit{symmetric} version.
\subsection{Non-symmetric HUP}
Observe that the HUP is a particular instance of the family of inequalities
\begin{equation}
\label{eq:uncerxiao}
    \norm{f}_q \lesssim \norm{f |x| ^{\alpha}}_q +\norm{\widehat{f} |\xi| ^\beta}_p
\end{equation}
for the parameters $p=q=2$ and $\alpha,\beta =1$.  Our first main goal is to fully characterize inequality
\eqref{eq:uncerxiao}. Before presenting the characterization of \eqref{eq:uncerxiao}, we would like to mention some facts about these inequalities.

This line of research was initiated by Cowling and Price~\cite{cowlingprice}, who studied an estimate similar to \eqref{eq:uncerxiao} and proved that,
 for $d=1$, $1\le p, q \le \infty$ 
 and $\alpha,\beta>  0$, the inequality
$$\norm{f}_2 \lesssim \norm{f |x| ^{\alpha}}_q +\norm{\widehat{f} |\xi| ^\beta}_p$$
holds if and only if $\alpha>\frac12-\frac1q$ and
$\beta>\frac12-\frac1p.$ For the multivariate version, see \cite{Luef, Tang}.

Note that  particular cases of \eqref{eq:uncerxiao}
were studied by  Steinerberger \cite{Steinerberger}
and 
Xiao \cite{xiao2022lp}, who
considered the cases $q=1$, $p=\infty$ and 
$p=q'$, $1 \leq q < \infty$, respectively. Both authors gave sufficient conditions for \eqref{eq:uncerxiao} to hold.

By homogeneity, it is easy to see that \eqref{eq:uncerxiao} is equivalent to 
\begin{equation}
\label{eq:uncerxiaoprim}
    \norm{f}_q\lesssim \norm{f |x| ^{\alpha}}^\frac{\beta+d(\frac1p+\frac1q-1)}{\alpha+\beta+d(\frac1p+\frac1q-1)}_q \norm{\widehat{f} |\xi| ^\beta}^\frac{\alpha}{\alpha+\beta+d(\frac1p+\frac1q-1)}_{p}
\end{equation} provided that $\beta\geq d(1-\frac{1}{q}- \frac{1}{p}).$
(In fact we will see that neither \eqref{eq:uncerxiao} nor \eqref{eq:uncerxiaoprim} is  valid for $\beta<d(1-\frac{1}{q}- \frac{1}{p}).$) 

We now proceed to state the characterization of \eqref{eq:uncerxiao}.
\begin{theorem}
\label{theo:mikxiao}
Let $1 \leq p,q \leq \infty$  and $\alpha,\beta>  0$.  Then,  inequality
\eqref{eq:uncerxiao}
holds if and only if one of the following conditions holds: 
\begin{itemize}
    \item $1\le p\leq q<\infty$ and $\beta=d(1-\frac{1}{q}- \frac{1}{p});$
    \item $\beta>d(1-\frac{1}{q}- \frac{1}{p}),$
\end{itemize}
and one of the following conditions holds:
\begin{itemize}
\item $q \geq 2;$
\item $p,q \leq 2;$
\item  $q<2 <p$\; and \; $\beta> d(\frac{1}{2}- \frac{1}{p})$.
\end{itemize}
\end{theorem}
\begin{remark}
   1. If 
    $q=1$ and  $p=\infty$, Theorem
    \ref{theo:mikxiao}
    recovers the sufficient condition given by Steinerberger. 
Besides, if    $p=q'$, $1 \leq q < \infty$, and $\alpha,\beta>  0$,
we see that 
\begin{equation}
\label{eq:nonsymp and p'}
\norm{f}_q\lesssim \norm{f |x| ^{\alpha}}^\frac{\beta}{\alpha+\beta}_q \norm{\widehat{f} |\xi| ^\beta}^\frac{\alpha}{\alpha+\beta}_{q'}, \quad
1 \leq q \leq \infty,
\end{equation}
holds if and only if
$$\beta>\begin{cases}
d(\frac{1}{q}-\frac{1}{2})\quad\mbox{if}\quad q<2,\\
0\qquad\qquad\,\mbox{if}\quad q\geq 2.
\end{cases}$$ This improves the sufficient  conditions by  
  Xiao in \cite{xiao2022lp} for  $q\geq 2$ as well as those   given by Price in \cite[Corollary 5.1]{price-studia} for $p=q=2$ and $\alpha=\beta$.

 2. Note that the classical 
 Hausdorff--Young
 estimate $\norm{f}_q\lesssim \norm{\widehat{f} }_{q'}$, $2 \leq q\leq\infty$, coincides with  
inequality \eqref{eq:nonsymp and p'}
 with $\beta=0$.
 In fact, Theorem
\ref{theo:mikxiao} allows us to consider the following weighted  Hausdorff--Young type
 estimate with
$\beta>0$:
\begin{equation}
\label{vspom0}
\norm{f}_q\lesssim 
\norm{\frac{f |x| ^{\alpha}}{\norm{f}_q} 
}_q ^\frac{\beta}{\alpha}
\norm{\widehat{f} |\xi| ^\beta}_{q'}, \quad
2 \leq q \leq \infty.
\end{equation}

 3. Another interesting case is when $p=q.$ The well-known Hardy--Littlewood estimate states that, for $q\ge 2,$
$$    \norm{f}_q \lesssim \norm{\widehat{f}
|\xi| ^\beta
}_q, \qquad 
\beta=
d(1-\frac2q),
$$
which is a particular case of \eqref{eq:uncerxiaoprim}.
Similarly as for the Hausdorff--Young inequality, Theorem
\ref{theo:mikxiao} allows us to obtain the following generalization of the Hardy--Littlewood inequality for $\beta>d(1-\frac2q)$   $$\norm{f}_q\lesssim \norm{f |x| ^{\alpha}}^\frac{\beta+
d(\frac2q-1)
}{\alpha+\beta+d(\frac2q-1)
    }_q \norm{\widehat{f} |\xi| ^\beta}^\frac{\alpha}{\alpha+\beta+d(\frac2q-1)
    }_{q},
    $$
    or, equivalently,
\begin{equation*}
\norm{f}_q\lesssim \norm{\frac{f |x| ^{\alpha}}{\norm{f}_q}}^\frac{\beta+
d(\frac2q-1)
}{\alpha
    }_q \norm{\widehat{f} |\xi| ^\beta}_{q},\qquad 1\le q\le \infty.
\end{equation*}

\end{remark}
In fact, in Section \ref{section:non} we will obtain Theorem \ref{theorem:GH},  a more general version of Theorem 
\ref{theo:mikxiao}, where we replace in inequality \eqref{eq:uncerxiao} the power weights $|x|^\alpha, |\xi|^\beta$ by broken power weights, which are defined by
\begin{equation}
\label{def:broken}
    x^A = \begin{cases}
|x|^{A_1}\quad\mbox{if}\quad |x|<1,\\
|x|^{A_2}\quad\mbox{if}\quad |x|\geq 1,\\
\end{cases}
\end{equation}
where $A=(A_1,A_2)$ with  $A_1,A_2\geq 0$. Throughout the paper we will use capital letters in the exponent
to denote broken power weights.

The reason for this generalization is that broken 
power weights allow us to clearly see 
the balance between the behavior of a function and its Fourier transform at the origin and infinity.

To illustrate our results, we give a simple corollary of Theorem \ref{theorem:GH} in the case $p=q'$.
\begin{corollary}\label{corollary-broken}
    Let $q\geq 2$ and $A_1,A_2,B_1,B_2>0$. Then, 
\begin{equation}
    \label{vspom1}
\norm{f}_q\lesssim 
\norm{\frac{f |x| ^{A}}{\norm{f}_q} 
}_q ^\frac{B^T}{A}
\norm{\widehat{f} |\xi| ^B}_{q'}, \quad
2 \leq q \leq \infty,
\end{equation} where $\frac{B^T}{A}=(\frac{B_2}{A_1},\frac{B_1}{A_2})$,
\end{corollary}

Observe that, in \eqref{vspom1}, $B_1$ and $B_2$ are related to $A_2$ and $A_1$, respectively. This is a manifestation of the interplay 
 between the behavior of $f$ and $\widehat{f}$ at zero and infinity, which cannot be seen in \eqref{vspom0}.

\subsection{Symmetric HUP}The symmetry between $f$ and $\widehat{f}$ is an appealing feature of the HUP, which is not present in \eqref{eq:uncerxiao}. 
Motivated by this, we consider the family of symmetric inequalities
\begin{equation}
\norm{f}_{p}\norm{\widehat{f}}_{p'} \lesssim \norm{|x|^\alpha f}_p \norm{|\xi|^\beta \widehat{f}}_{p'},
    \label{ineq:uncer}
\end{equation} to which the classical Heisenberg principle clearly belongs.
We will see that  
\eqref{ineq:uncer} holds
  if and only if  $\alpha=\beta>d(\frac{1}{2}- \frac{1}{p})$.

 Once again, in order to highlight the aforementioned 
 interplay 
 between $f$ and $\widehat{f},$
 we  will replace the power weights by
broken power weights. This leads us to the following result,
 \begin{theorem}
 \label{theorem:broken}
Let $2\le p\le \infty$ and $A_1,A_2,B_1,B_2> 0$. Then, the
inequality \begin{equation}   \norm{f}_{p}\norm{\widehat{f}}_{p'} \lesssim \norm{|x|^A f}_p \norm{|\xi|^B \widehat{f}}_{p'}.
    \label{ineq:uncer--}
\end{equation} holds if and only if $B_2 \geq A_1, A_2 \geq B_1$, $A_2>d(\frac{1}{2}-\frac{1}{p})$ and one the following conditions holds:
\begin{enumerate}
    
    \item $A_1<d(\frac{1}{2}-\frac{1}{p})$ and $B_2 ^2\geq d(\frac{1}{2}-\frac{1}{p})A_1;$
    \item $A_1>  d(\frac{1}{2}-\frac{1}{p});$
    \item $ A_1=d(\frac{1}{2}-\frac{1}{p})  <B_2.$
  
\end{enumerate}

 \end{theorem}

Note that replacing $f$ by $\widehat{f}$ in inequality \eqref{ineq:uncer--} has the effect of exchanging the roles of $p$ by $p'$ and $A$ by $B$, so the condition $p\geq 2$ is not restrictive.

Finally,  we would like to mention 
the recent  
 approach suggested by Wigdersons in \cite{wig} and elaborated in \cite{dias,Tang}, which establishes 
Heisenberg-type uncertainty principles for general operators sharing some properties with the Fourier transform.  
In Section \ref{section: special} we will discuss how our results fit in this framework.

\subsection{Broken power weights}
Here we collect some simple properties of broken power weights.

For $A=(A_1,A_2)$ and
$B=(B_1,B_2)$,
we define 
$A\pm B$, $AB$, 
$A/B$, and $A\le B$ coordinate-wise. 
Observe that $(1/x)^A=x^{-A^T}$   and
$(x^B)^A=x^{B A}$, in particular, 
$(x^\alpha)^A=x^{\alpha A}$, $\alpha\in \mathbb{R}$, but
in general $(xy)^A\neq x^Ay^A$.

We will also use that
\begin{equation}
\label{eq:integral}\int_0 ^x y^A dy\approx \begin{cases}
    x^{A+1}, &A_1,A_2>-1,\\
    x^{A+1} + \log(x+1), &A_1>A_2=-1,\\
    x^{A+1} + \mathbbm{1}_{x>1}, &A_2<-1<A_1,\\
    \infty, &A_1\leq -1,
\end{cases}\end{equation}
and
\begin{equation}
\label{eq:integral2}\int_x^\infty y^{-A} dy= \begin{cases}
    x^{-A+1} +\mathbbm{1}_{x<1}, &A_1<1<A_2,\\
    x^{-A+1} + \log_+(1/x), &A_1=1<A_2,\\
    x^{-A+1} , &A_1,A_2>1,\\
    \infty, &A_2\leq 1.
\end{cases}\end{equation}
\subsection{Structure of the paper}
We now explain the organization of the remainder of this paper.
In Section \ref{section:local}, we characterize 
  the optimal behavior of $C_{q,p}(t)$  in the inequality
  $\displaystyle{\Big(\int_{|\xi|\leq t} |f|^q\Big)^{\frac{1}{q}}}\le C_{q,p}(t){ \norm{\widehat{f}|\xi|^D}_p}
$
and study a similar problem for trigonometric polynomials. This is a well-known problem in harmonic analysis, referred to as a local uncertainty estimate.

In Section 
\ref{section:non}, we study the non-symmetric HUP mentioned in the introduction. More precisely, we find the optimal function  $H(\cdot)$ so that $\norm{f}_q \lesssim H(x) $
whenever
$   \norm{\widehat{f}|\xi|^\beta}_p\leq 1$  and 
$\norm{{f}|t|^\alpha}_q\leq x$, { for } $x\in \mathbb{R}^+.$
 This result implies Theorem
\ref{theo:mikxiao} and Corollary \ref{corollary-broken}.
 Section \ref{section:broken}
 contains the characterization of the symmetric HUP, Theorem
\ref{theorem:broken}.

In Section \ref{section: special}
 we discuss our results in Wigderson's framework.
Finally, in the Appendix we state the necessary and sufficient conditions  for the weighted Fourier inequality
$\big\|u\widehat{f}\,\big\|
_{q}\lesssim\norm{ f v}_{p}
$
to be true.

\section{Local  uncertainty principle}\label{section:local}
The study  of local Fourier inequalities  $$    \left(\int_{|\xi|\leq t} |f|^2\right)^{\frac{1}{2}}\lesssim  C(t) \norm{\widehat{f}|\xi|^\beta}_p, \qquad \beta>0,$$ was initiated by Faris
\cite{faris} and
Price  \cite{price}, and further developed in \cite{ricci,price-studia}, among many other works (see also \cite[Section 4]{folland}). We note that this inequality is a particular case of a weighted Fourier inequality with monotonic weights, whose complete characterization has been recently been obtained by the authors in \cite{saucedo2025fouriertransformextremizerclass}.

 The lemma below contains the solution of 
  the question posed by Price 
 both for any $1\leq p,q\leq \infty$ and  broken power weights. For given $p$ and $q$, we define $\frac{1}{r}=\frac{1}{q}-\frac{1}{p}$.
 \begin{lemma}
\label{lemma:local} Let $D=(D_1,D_2)\geq 0$ and define

\begin{equation}
\label{ineq:localD2}
    C_{q,p}(t)=\sup_{f} \frac{\left(\int_{|\xi|\leq t} |f|^q\right)^{\frac{1}{q}}}{ \norm{\widehat{f}|\xi|^D}_p},
\end{equation}

Then, 
\begin{enumerate}

    \item Case $1<p\leq q<\infty$. If $D_2\geq d(1-\frac{1}{p}-\frac{1}{q})$ and $D_1<\frac{d}{p'}$,

\begin{align*}&\mbox{for }t<1, \; C_{q,p}(t)\approx\begin{cases}
      t^{D_2+ d(\frac{1}{p}+\frac{1}{q}-1)},& D_2<\frac{d}{p'},\\
  t^{\frac{d}{q}} \log^{\frac{1}{p'}}_+(1/t), & D_2= \frac{d}{p'},
   \\
    t^{\frac{d}{q}}, &D_2>\frac{d}{p'};\\    
    \end{cases}\\
&\mbox{for }t>1, \; C_{q,p}(t)\approx\begin{cases}
      t^{D_1+ d(\frac{1}{p}+\frac{1}{q}-1)},& D_1\geq d(1- \frac{1}{p}-\frac{1}{q}),\\
    1, &D_1<d(1- \frac{1}{p}-\frac{1}{q}).\\   
    \end{cases}
    \end{align*}
        \item Case  $1\leq q<p\leq \infty,$ $ \max(q,p')\geq 2$. If $D_2> d(1-\frac{1}{p}-\frac{1}{q})$ and $D_1<\frac{d}{p'}$,    \begin{align*}   
&\mbox{for }t<1, \;C_{q,p}(t)\approx\begin{cases}
      t^{D_2+ d(\frac{1}{p}+\frac{1}{q}-1)}, & D_2<\frac{d}{p'},\\
  t^{\frac{d}{q}} \log^{\frac{1}{p'}}_+(1/t), & D_2= \frac{d}{p'},
   \\
    t^{\frac{d}{q}}, &D_2>\frac{d}{p'};\\    
    \end{cases}\\            
&\mbox{for }t>1, \;C_{q,p}(t)\approx\begin{cases}
      t^{D_1+ d(\frac{1}{p}+\frac{1}{q}-1)},  &D_1>d(1- \frac{1}{p}-\frac{1}{q}),\\
    \log^{\frac{1}{r}}_+(t), & D_1=d(1- \frac{1}{p}-\frac{1}{q}),
   \\
    1, &D_1<d(1- \frac{1}{p}-\frac{1}{q}).\\
    \end{cases}
    \end{align*}

   \item Case  $1\leq q<2<p \leq \infty$. If $D_2> d(\frac{1}{2}-\frac{1}{p})$ and $D_1<\frac{d}{p'}$,        
\begin{align*}
&\mbox{for }t<1, \;C_{q,p}(t)\approx \begin{cases}
      t^{D_2+ d(\frac{1}{p}+\frac{1}{q}-1)},  &D_2<\frac{d}{p'},\\
   t^{\frac{d}{q}} \log^{\frac{1}{p'}}_+(1/t), & D_2=\frac{d}{p'},
   \\
    t^{\frac{d}{q}}, &D_2>\frac{d}{p'};    
    \end{cases}\\ 
&\mbox{for }t>1, \;C_{q,p}(t)\approx \begin{cases}
      t^{D_1+ d(\frac{1}{p}+\frac{1}{q}-1)},  &D_1>d( \frac{1}{2}-\frac{1}{p}),\\
    t^{d(\frac{1}{q}-\frac{1}{2})}\log^{\frac{1}{2}-\frac{1}{p}}_+(t), & D_1=d( \frac{1}{2}-\frac{1}{p}),
   \\
    t^{d(\frac{1}{q}-\frac{1}{2})}, &D_1<d(\frac{1}{2}-\frac{1}{p}).
    \end{cases}\end{align*} 
   \item Case $p=1$ and $D_1=0$. $\,\,C_{q,p}(t)\approx t^{\frac{d}{q}} $.
\item Case $q=\infty$ and $p>1$.
If $D_1<\frac{d}{p'}<D_2$,
$C_{q,p}(t)\approx 1$.
    \end{enumerate}
In all cases not covered above, $C_{q,p}(t)=\infty$.
\end{lemma}
\begin{proof}

By Theorem \ref{theorem:mainh}, in the first, fourth and fifth cases,
\begin{align*}
    C_{q,p}(t)&\approx  \sup_{0\leq s \leq t} \left(\int_{0}^s x^{d-1} dx \right)^{\frac{1}{q}} \left(\int_{0}^{\frac{1}{s}} x^{-p'D} x^{d-1} dx \right)^{\frac{1}{p'}}.
\end{align*}
For the fourth and fifth case, the result follows easily. In the first case, a computation shows that
\begin{align*}
  C_{q,p}(t)
    \approx \sup_{0\leq s \leq t} \begin{cases}
        s^{d(\frac{1}{q}+\frac{1}{p}-1)+D^T}, &D_1,D_2<\frac{d}{p'},\\
        s^{d(\frac{1}{q}+\frac{1}{p}-1)+D^T} + s^{\frac{d}{q}} \log^{\frac{1}{p'}}_+(1/s), &D_1<\frac{d}{p'}=D_2,\\
        s^{d(\frac{1}{q}+\frac{1}{p}-1)+D^T} + s^{\frac{d}{q}} \mathbbm{1}_{s<1}, &D_1<\frac{d}{p'}<D_2,\\
        \infty, &D_1\geq \frac{d}{p'}.
    \end{cases}
\end{align*}

In the second case, Theorem \ref{theorem:mainh} implies that
\begin{align*}
    C_{q,p}(t)&\approx Q_{q,p}(t):=\left(\int_{0}^t s ^{d-1} \left(\int_{0}^s x^{d-1} dx \right)^{\frac{r}{p}} \left(\int_{0}^{\frac{1}{s}} x^{-p'D} x^{d-1} dx \right)^{\frac{r}{p'}} \right)^{\frac{1}{r}}\\
    &\approx \begin{cases}
        \left(\int_{0}^t \big(s^{d(\frac{1}{q}+\frac{1}{p}-1)+D^T}\big)^r \frac{ds}{s}\right)^{\frac{1}{r}}, &D_1,D_2<\frac{d}{p'},\\
        \left( \int_0 ^t \big(s^{d(\frac{1}{q}+\frac{1}{p}-1)+D^T} + s^{\frac{d}{q}} \log^{\frac{1}{p'}}_+(1/s)\big)^r {\frac{ds}{s}}\right)^{\frac{1}{r}}, &D_1<\frac{d}{p'}=D_2,\\
        \left(\int_0 ^t   \big(s^{d(\frac{1}{q}+\frac{1}{p}-1)+D^T} + s^{\frac{d}{q}} \mathbbm{1}_{s<1}\big)^r {\frac{ds}{s}}\right)^{\frac{1}{r}}, &D_1<\frac{d}{p'}<D_2,\\
        \infty, &D_1\geq \frac{d}{p'}.
    \end{cases}
\end{align*}

In the third case, instead of using Theorem \ref{theorem:mainh} we obtain the upper bound by using Hölder's inequality and the second item as follows:
$$\left(\int_{|\xi|\leq t} |f|^q\right)^{\frac{1}{q}}\lesssim t^{d \left(\frac{1}{q}- \frac{1}{2}\right)}\left(\int_{|\xi|\leq t} |f|^2\right)^{\frac{1}{2}}\lesssim  t^{d \left(\frac{1}{q}- \frac{1}{2}\right)} C_{2,p}(t) \norm{\widehat{f}|\xi|^D}_p,$$ 
For the lower bound, in view of Theorem \ref{theorem:mainh}, we first note that
$Q_{p,q}(t)\lesssim C_{p,q}(t)$, which gives the result for $D_1\geq \frac{d}{p'}$ and for $t<1$, $D_2>d(\frac{1}{2}-\frac{1}{p})$. Second, by Lemma \ref{lemma:prelinec}, for all $y>0$,
$$y^{-\frac{d}{2}} \left(\int_{0}^y \mathbbm{1}_{[0,t]}(x) x^{d-1} dx \right)^{\frac1q} \left(\int_{1/y}^\infty x^{-p^\# D } x^{d-1} dx \right)^{\frac12 -\frac1p} \lesssim C_{q,p}(t). $$ Setting $y=t/2$ we deduce that
\begin{equation}
C_{q,p}(t) \gtrsim \begin{cases}
   t^{D^T + d (\frac1p+\frac1q-1)}, &D_1,D_2>d(\frac{1}{2}-\frac{1}{p}),\\
   t^{D^T + d (\frac1p+\frac1q-1)} + t^{d(\frac{1}{q}-\frac12)} \log^{\frac{1}{2}-\frac{1}{p}}_+(t), &D_2>d(\frac{1}{2}-\frac{1}{p})=D_1,
   \\
    t^{D^T + d (\frac1p+\frac1q-1)} + t^{d(\frac{1}{q}-\frac12)} \mathbbm{1}_{t>1},&D_2>d(\frac{1}{2}-\frac{1}{p})>D_1,\\
    
    \infty, &\mbox{ otherwise}
    ,
\end{cases}\end{equation} whence the result follows in the remaining cases. 

\end{proof}

We will also need the local Fourier estimates for series. For 
$T_M(x)=\sum_{|n|\leq M} c_n e^{2 \pi i \langle x, n\rangle}$, we are interested in the optimal $C_{q,p}(M)$ so that 
\begin{equation}\label{local}
\|T_M\|_q\le C_{q,p}(M){\left(\sum_{|n|\leq M} |c_n|^p(1+|n|)^{p\gamma}\right)^{\frac{1}{p}}}.
\end{equation}The classical Pitt inequality \cite{pitt, stein} ensures that 
$C_{q,p}(M)\lesssim1$ for some range of the parameters.
Price and Racki
\cite{Price-Racki} characterized the dual  to \eqref{local} when $p=2$.
Later, Foschi \cite{foschi}
found the optimal behavior of $C_{q,p}(M)$ for the case $\gamma=0$. Note also that estimate \eqref{local}
 for $p=2$ is closely related to the Nikol'skii (and Bernstein-Nikol'skii) inequality:
$$\|T_M\|_q\lesssim M^{d(\frac1p-\frac1q)}\|T_M\|_p, \qquad 0<p\le q\le\infty,$$  see, e.g.,  \cite{jfaa}.
In the next result we characterize \eqref{local} for all parameters.
\begin{lemma}
\label{lemma:X} Let $\gamma\geq 0$ and define
$$C_{q,p}(M)= \sup_{c_n} \frac{ \left(\int_{[0,1]^d}\left | \sum_{|n|\leq M} c_n e^{2 \pi i \langle x, n\rangle}\right|^q\right)^{\frac{1}{q}}}{\left(\sum_{|n|\leq M} |c_n|^p(1+|n|)^{p\gamma}\right)^{\frac{1}{p}}}.$$
Then,
\begin{enumerate}
    \item for $1\le p\leq q<\infty$, 
$$C_{q,p}(M) \approx 
\begin{cases}
M^{-\gamma +d(1 -\frac{1}{p}- \frac{1}{q})}, &\gamma < d(1- \frac{1}{p}- \frac{1}{q}),\\
1, &\gamma \geq d(1- \frac{1}{p}- \frac{1}{q});
\end{cases}$$
\item for 
$1\le p\leq q=\infty$,
$$C_{\infty,p}(M) \approx \begin{cases}
M^{-\gamma +d(1 -\frac{1}{p})}, &\gamma < d(1- \frac{1}{p}),\\
\log(M+1)^{1-\frac{1}{p}}, &\gamma = d(1- \frac{1}{p}), \\
1, &\gamma > d(1- \frac{1}{p});
\end{cases}$$
\item for $2\leq q<p \leq \infty$, 

$$C_{q,p}(M)\approx
\begin{cases}
M^{-\gamma +d(1 -\frac{1}{p}- \frac{1}{q})}, &\gamma < d(1- \frac{1}{p}- \frac{1}{q}),\\
\log(M+1)^{\frac{1}{q}-\frac{1}{p}}, &\gamma = d(1- \frac{1}{p}- \frac{1}{q} ),\\
1, &\gamma > d(1- \frac{1}{p}- \frac{1}{q});
\end{cases}$$

\item for $1\leq q  < 2$, \,\,$C_{q,p}\approx C_{2,p}.$
\end{enumerate}
\end{lemma}
\begin{proof}

Observe that defining $1/v_n=(1+|n|)^ {-\gamma} \mathbbm{1}_{|n|\leq M}$, we have
$$C_{q,p}(M)= \sup _c \frac{\norm{f}_q}{\norm{\widehat{f}v}_p}.$$ Thus, the version of Theorem \ref{theorem:mainh} for Fourier series (see Remark 5.4 in \cite{saucedo2025fouriertransformextremizerclass}) shows that, for $1\leq p\leq q \leq  \infty$,
$$C_{q,p}(M) \approx \sup_{1\leq n\leq M} \left(\int_{|x|\leq n ^{-1}} dx \right)^{\frac{1}{q}} \left(\sum_{|j|\leq n} (1+|j|)^{-p'\gamma } \right)^{\frac{1}{p'}}.$$
A computation shows that, for $1\le p\le q<\infty$, 
\begin{align*}
 C_{q,p}(M)&\approx \sup_{1\leq n\leq M} \left(\int_{|x|\leq n ^{-1}} dx \right)^{\frac{1}{q}} \left(\sum_{|j|\leq n} (1+|j|)^{-p'\gamma } \right)^{\frac{1}{p'}}\\
 &\approx
\begin{cases}
M^{-\gamma +d(1 -\frac{1}{p}- \frac{1}{q})}, &\gamma < d(1- \frac{1}{p}- \frac{1}{q}),\\
1, &\gamma \geq d(1- \frac{1}{p}- \frac{1}{q}).
\end{cases}\end{align*}
Similarly, for $q=\infty$,
$$C_{\infty,p}(M)\approx \left(\sum_{|n|\leq M} (1+|n|)^{-p'\gamma} \right)^{\frac{1}{p'}} \approx \begin{cases}
M^{-\gamma +d(1 -\frac{1}{p})}, &\gamma < d(1- \frac{1}{p}),\\
\log(M+1)^{1-\frac{1}{p}}, &\gamma = d(1- \frac{1}{p}), \\
1, &\gamma > d(1- \frac{1}{p}).
\end{cases}$$

For $2\leq q<p\leq \infty$, Theorem \ref{theorem:mainh} shows that (recall  that $\frac{1}{r}=\frac{1}{q}-\frac{1}{p}$)
\begin{align*}
    C_{q,p}(M)&\approx \left(\sum_{|n|\leq M} (1+|n|)^{-p'\gamma } \left(\sum_{|j|\leq n} (1+|j|)^{-p'\gamma } \right)^{\frac{r}{q'}} \left(\int_{|x|\leq n^{-1}}  dx \right)^{\frac{r}{q}} \right)^{\frac{1}{r}}\\
    &\approx
\begin{cases}
M^{-\gamma +d(1 -\frac{1}{p}- \frac{1}{q})}, &\gamma < d(1- \frac{1}{p}- \frac{1}{q}),\\
\log(M+1)^{\frac{1}{r}}, &\gamma = d(1- \frac{1}{p}- \frac{1}{q} ),\\
1, &\gamma > d(1- \frac{1}{p}- \frac{1}{q}).
\end{cases}
\end{align*}
For $1\leq q \leq2$, by Hölder's inequality, $C_{q,p}(M)\leq C_{2,p}(M)$. For the converse, observe that Khinchin's inequality implies that for any sequence $c_n$, the average over all possible sign configurations $\varepsilon_n=\pm1$ satisfies $$\mathbbm{E}\left(\int_{[0,1]^d}\left | \sum_{|n|\leq M} \varepsilon_n c_n e^{2 \pi i \langle x, n\rangle}\right|^q\right) \approx \left(\int_{[0,1]^d}\left | \sum_{|n|\leq M}  c_n e^{2 \pi i \langle x, n\rangle}\right|^2\right)^{\frac q2}, $$ which implies that 
$C_{q,p}(M)\gtrsim C_{2,p}(M).$
\end{proof}

\section{Non-symmetric uncertainty inequalities
}\label{section:non}

\subsection{Main result}
Our goal in this section is to   extend \eqref{eq:uncerxiao}
 to the case of broken power weights.

Using the properties of $\norm{\widehat{f}|\xi|^\beta}_p$ and  $\norm{{f}|x|^\alpha}_q$ under the transformation $f(x)\mapsto Cf(\lambda x),$ where $\lambda ,C \in \mathbb{R}_+$, it is easy to see that \eqref{eq:uncerxiaoprim} is equivalent to
\begin{equation}
\label{eq:equivuncer}
    \norm{f}_q \lesssim 1 \mbox{ whenever } \norm{\widehat{f}|\xi|^\beta}_p\leq 1 \mbox{ and } \norm{{f}|x|^\alpha}_q\leq 1.
\end{equation}

 In the case of broken power weights, the dilation  
  is no longer nicely behaved, so the correct generalization of \eqref{eq:equivuncer} is
 \begin{equation}
\label{eq:equivuncer2}
    \norm{f}_q \lesssim H(x) \mbox{ whenever } \norm{\widehat{f}|\xi|^B}_p\leq 1 \mbox{ and } \norm{{f}|x|^A}_q\leq x, \mbox{ for } x\in \mathbb{R}^+.
\end{equation} 
Observe that, by considering translations and modulations,  
\eqref{eq:equivuncer2} is equivalent to 
\begin{equation}
\label{eq:equivuncer23}
    \norm{f}_q \lesssim H(x) \mbox{ whenever } \norm{\widehat{f}|\xi|^B}_p= 1 \mbox{ and } \norm{{f}|x|^A}_q= x, \mbox{ for } x\in \mathbb{R}^+.
\end{equation} 
Moreover, this can be rewritten 
 as
$$
{\norm{f}_q}\lesssim
\norm{\widehat{f} |\xi|^B}_p
H\Bigg(
\frac{\norm{f |x|^A}_q}{\norm{\widehat{f} |\xi|^B}_p}  \Bigg).
$$ 

We now  state the  characterization of the optimal function $H$.
\begin{theorem}
\label{theorem:GH}

Let $1 \leq p,q \leq \infty$, $A=(A_1,A_2)$ and $B=(B_1,B_2)$, with $A_1,A_2,B_1,B_2>0$. Define 
 $$H(x):=H_{A,B}(x)=\sup _{f \neq 0, \norm{f |x|^A}_q\leq x, \norm{\widehat{f} |\xi|^B}_p\leq 1}  {\norm{f}_q}$$ and
 $$\beta_{p,q}=\left(\max\Big(B_1+ d(\frac{1}{p}+\frac{1}{q}-1),0\Big),B_2+ d\big(\frac{1}{p}+\frac{1}{q}-1\big)\right).$$ Then,
\begin{enumerate}
    \item if 
    $1\leq p\le q<\infty$
      and $\frac{B_2}{d} \geq 1- \frac{1}{p}-\frac{1}{q}$,
    $$H(x)\approx x^{\frac{\beta_{p,q}^T}{A+\beta_{p,q}^T}};$$ 

    \item if $1\leq q<p\leq \infty$, $\max(q,p')\geq 2$, and $\frac{B_2}{d} > 1- \frac{1}{p}-\frac{1}{q}$,
    \begin{enumerate}
        \item if $B_1+ d(\frac{1}{q}+ \frac{1}{p}-1)\neq 0,$
        $$H(x)\approx x^{\frac{\beta_{p,q}^T}{A+\beta_{p,q}^T}};$$ 
        \item if $B_1+ d(\frac{1}{q}+ \frac{1}{p}-1)= 0,$
        $$H(x)\approx\begin{cases}
        x^{\frac{\beta_{p,q}^T}{A+\beta_{p,q}^T}}, &x<1, \\
       {\log(x+1)^{\frac{1}{r}}}, &x\geq 1.
        \end{cases}$$
    \end{enumerate}
  
   \item if $1\leq q<p\leq \infty$, $\max(q,p')<2$ and $\frac{B_2}{d} > 1- \frac{1}{p}-\frac{1}{2}$,
      \begin{enumerate}
        \item if $B_1+ d(\frac{1}{2}+ \frac{1}{p}-1)\neq 0,$
        $$H(x)\approx x^{\frac{\beta_{p,2}^T +d(\frac{1}{q} -\frac{1}{2})}{A+\beta_{p,2}^T +d(\frac{1}{q} -\frac{1}{2})}};$$ 
        \item if $B_1+ d(\frac{1}{2}+ \frac{1}{p}-1)= 0,$
        $$H(x)\approx \begin{cases}
        x^{\frac{\beta_{p,2}^T +d(\frac{1}{q} -\frac{1}{2})}{A+\beta_{p,2}^T +d(\frac{1}{q} -\frac{1}{2})}}, &x<1, \\
      \log(x+1)^{\frac{A_2 (\frac{1}{2}-\frac{1}{p})}{A_2+ d(\frac{1}{q}-\frac{1}{2})}}  x^{\frac{d(\frac{1}{q}-\frac{1}{2})}{A_2+ d(\frac{1}{q}-\frac{1}{2})}}, &x\geq 1.
        \end{cases}$$
    \end{enumerate}
    \item if $q=\infty$ and $B_2 > \frac{d}{p'}$,
   \begin{enumerate}
\item if $B_1+ d(\frac{1}{p}-1)\neq 0$, 
    
    $$H(x)\approx x^{\frac{\beta^T_{p,\infty}}{A+\beta^T_{p,\infty}}},$$
    
    \item if $B_1+ d(\frac{1}{p}-1)=0$, 
    $$H(x)\approx \begin{cases}
        x^{\frac{\beta^T_{p,\infty}}{A+\beta_{p,\infty}^T}}, &x<1, \\
        {\log(x+1)^{\frac{1}{p'}}}, &x\geq 1.
        \end{cases}$$
\end{enumerate}
    \item  In all cases not covered above, $H(x)=\infty$.
    \end{enumerate}
   
\end{theorem}

 Theorem \ref{theo:mikxiao} follows from Theorem \ref{theorem:GH} by noting that in the case of non-broken weights either $H(x)=x^{\frac{\beta + d (\frac1p+\frac1q-1)}{\alpha +\beta + d (\frac1p+\frac1q-1)}}$ or $H=\infty$.

The proof of Theorem \ref{theorem:GH} will be given in the next two subsections.

\subsection{Proof of $\lesssim$ part in Theorem \ref{theorem:GH}}

First, we prove Lemmas \ref{lemma:local2} and \ref{lemma:interpq}, which relate the size of $\norm{\widehat{f}|\xi|^\beta}_p$ with the concentration of $f$ around a point.

\begin{lemma}
\label{lemma:local2}
Let $1\leq p,q\leq \infty,$ $ q<\infty,$ and $A,B\geq0$. Then, for $0\leq D\leq B$, 
\begin{equation}
\label{ineq:localdiff}
\norm{f}_q \lesssim C(t) t^{L^T} \norm{\widehat{f}|\xi|^B}_p + t^{-A} \norm{f |x|^A}_q, \qquad t>0,
\end{equation}
where  $L=B-D\geq 0$, and $C:=C_{q,p}$ is defined in \eqref{ineq:localD2}.
\end{lemma}
\begin{proof}
Before starting the proof, we define $$\Delta_t f(x)=f(x)-f(x+te_1),\qquad \Delta_t^{k} f=
\Delta_t (\Delta_t^{k-1} f),\quad k\in\mathbb{N},$$ and note that
$$(\Delta_h^Nf )^{\wedge}(\xi) = \widehat{f}(\xi) (1-e^{2 \pi i  \xi_1 h})^N,\quad N\in\mathbb{N}.$$ First, recalling definition \eqref{ineq:localD2} in Lemma \ref{lemma:local}, we see that
$$\left(\int_{|\xi|\leq t} |\Delta^N_{2t}f|^q\right)^{\frac{1}{q}} \lesssim C_{q,p}(t) \norm{\widehat{f}(\xi)|\xi|^D (1-e^{2 \pi i  \xi_1 t})^N}_{p}.$$ Next, we use 
$$\left |(1-e^{2 \pi i  \xi_1 t})^N \right | \lesssim \min(1, t| \xi|)^N\leq t^{L^T}  |\xi|^L$$ provided $N\geq 
\max(L_1, L_2)
$
to conclude that
$$\left(\int_{|\xi|\leq t} |\Delta^N_{2t}f|^q\right)^{\frac{1}{q}}   \lesssim C_{q,p}(t)t^{L^{T}} \norm{
\widehat{f}(\xi)|\xi|^{D+L}}_p.$$

Second, the triangle inequality gives
$$\left(\int_{|x|\leq t} |f|^q \right)^{\frac{1}{q}} \lesssim \left(\int_{|x|\leq t} | \Delta^N_{2t}f |^q \right)^{\frac{1}{q}}+\left(\int_{|x|\leq t} |f-\Delta^N_{2t}f  |^q \right)^{\frac{1}{q}}.$$
Third, we observe that
\begin{align*}
    \int_{|x|\leq t} |f-\Delta^N_{2t}f  |^q &=\int_{|x|\leq t} \left | \sum_{j=1}^N (-1)^j \binom{N}{j} f(x+2jte_1)  \right|^q\\
&\lesssim \sum_{j=1}^N  \int_{|x|\leq t} |f(x+2jte_1)|^q dx \\
&\lesssim \int_{|x|\geq t} |f|^q.
\end{align*} 
Finally, the result follows by noting that
\begin{equation*}
     \int_{|x|\geq t} |f|^q \lesssim t^{-Aq}\norm{f |x|^A}_q^q. 
\end{equation*}
\end{proof}
The counterpart of the previous lemma for  $q=\infty$ reads as follows.
\begin{lemma}
\label{lemma:interpq}
Let $1\leq p\leq \infty$ and $A,B\geq 0$. Assume also that if $p>1$,  $p'B_2>d$. Then \begin{equation}
\label{eq:infinity}
    \norm{f}_\infty \lesssim C(t) \norm{\widehat{f}|\xi|^{B}}_{p} +t^{-A} \norm{f |x|^A}_\infty,
\end{equation}
where 
\begin{enumerate}
    \item if $B_1 + d(\frac{1}{p}-1) \neq 0,$
    $$C(t)=t^{\beta_{p,\infty}^T},$$
    \item if $B_1+ d(\frac{1}{p}-1)=0$,
    $$C(t)= \begin{cases}
    t^{{\beta_{p,\infty}^T}}, &t<1,\\
    \log(1+t)^{\frac{1}{p'}}, &t\geq 1.
    \end{cases}$$ 
\end{enumerate}
\end{lemma}

\begin{proof}
Let $\Phi$ be a Schwartz function such that $\widehat{\Phi}$ is a compactly supported function which equals  $1$ on the unit ball and $0$ outside of the ball of radius 2.

We write
\begin{equation*}
    \begin{aligned}
    f(x)&=
    \int_{\mathbb{R}^d} e^{2 \pi i \langle \xi,y \rangle} \widehat{f}(\xi) d \xi\\
    &= \int_{\mathbb{R}^d} e^{2 \pi i \langle \xi,y \rangle} \widehat{\Phi}(\xi/h) \widehat{f}(\xi) d \xi+ \int_{\mathbb{R}^d} e^{2 \pi i \langle \xi,y \rangle} (1-\widehat{\Phi}(\xi/h)) \widehat{f}(\xi) d \xi\\
    &=J_1+J_2.
    \end{aligned}
\end{equation*}

For $J_2$, since $\widehat{\Phi}(\xi)=1$ for $|\xi|\leq 1$, by using Hölder's inequality we deduce that
\begin{equation*}
    \begin{aligned}
       J_2&= \int_{|\xi|\geq t} e^{2 \pi i \langle \xi,y \rangle} (1-\widehat{\Phi}(\xi/t)) \widehat{f}(\xi) d \xi \lesssim \int_{|\xi|\geq t} |\widehat{f}(\xi)| d \xi\\
        &\lesssim\left(\int_{|\xi|\geq t} |\xi|^{-p'B}d \xi\right)^{\frac{1}{p'}} \norm{\widehat{f}|\xi|^{B}}_{p}.
    \end{aligned}
\end{equation*}

For $J_1$, we first note that, for $N\geq A$, 
\begin{equation}\label{eq:qqq}
    \begin{aligned}
        \left |\int_{\mathbb{R}^d} \widehat{\Phi}(\xi/t) e^{-2 \pi i \langle y , \xi \rangle} \Delta^N_{3t}  \widehat{f}(\xi)d \xi\right|&= \left |t^d \int_{\mathbb{R}^d} \Phi(t(x-y)) (1-e^{2 \pi i \langle 3t, x\rangle})^N  f(x) dx \right|\\
        & \lesssim t^{d+A^T} \norm{f |x|^A}_\infty \int_{\mathbb{R}^d} |\Phi(t(x-y))|  dx \\
        &\approx t^{A^T} \norm{f |x|^A}_\infty,
    \end{aligned}
\end{equation}
where we used the estimate
$$ \left|1-e^{2 \pi i \langle 3t, x \rangle} \right|^N\lesssim  t^{A^T} x^A.$$

Second, since $\widehat{\Phi}$ is supported in the ball of radius $2$, using inequality \eqref{eq:qqq}, we arrive at 
\begin{equation*}
    \begin{aligned}
        J_1
&= \int_{|\xi|\leq 2t} e^{2 \pi i \langle \xi,y \rangle} \widehat{\Phi}(\xi/t) \left(\Delta^N _{3t} \widehat{f}(\xi)-\Delta^N _{3t} \widehat{f}(\xi) + \widehat{f} (\xi)\right) d \xi \\
&\lesssim \int_{|\xi|\geq t} |\widehat{f}(\xi)| d \xi + t^{A^T} \norm{f |x|^A}_\infty\\
&\lesssim \left(\int_{|\xi|\geq t} |\xi|^{-p'B}d \xi\right)^{\frac{1}{p'}} \norm{\widehat{f}|\xi|^{B}}_{p}+ t^{A^T} \norm{f |x|^A}_\infty.
    \end{aligned}
\end{equation*}
In conclusion, 
$$\left |\int_{\mathbb{R}^d} e^{2 \pi i \langle \xi,y \rangle} \widehat{f}(\xi) d \xi \right| \lesssim   \left(\int_{|\xi|\geq t} |\xi|^{-p'B}d \xi\right)^{\frac{1}{p'}} \norm{\widehat{f}|\xi|^{B}}_{p}+ t^{A^T} \norm{f |x|^A}_\infty.$$
Equivalently,
$$\norm{f}_\infty \lesssim C(t) \norm{\widehat{f}|\xi|^{B}}_{p} +t^{-A} \norm{f |x|^A}_\infty,$$
where 
\begin{enumerate}
    \item if $B_1 + d(\frac{1}{p}-1) \neq 0$,\quad
$C(t)=t^{\beta_{p,\infty}^T},$
    \item if $B_1+ d(\frac{1}{p}-1)=0,$\quad
    $C(t)= \begin{cases}
    t^{{\beta_{p,\infty}^T}}, &t<1,\\
    \log(1+t)^{\frac{1}{p'}}, &t\geq 1.
    \end{cases}$
\end{enumerate}
\end{proof}

In order to complete the proof of sufficiency,  all that remains is to choose the right values of $D$ and $t$ in Lemmas \ref{lemma:local2} and \ref{lemma:interpq}.
\begin{proof}[Proof of $\lesssim$ part in Theorem \ref{theorem:GH}]\,
By Lemmas \ref{lemma:local2}, \ref{lemma:interpq} and the condition  $\norm{\widehat{f}|\xi|^B}_p\le1$, $$\norm{f}_q\lesssim C(t) t^{L^T} + t^{-A}\norm{fx^A}_q, \qquad t>0,$$ where $D\leq B$ and $L=B-D$. We now proceed to choose convenient values of $D$ and $t$ in each case.

\begin{enumerate}

    \item Let  $1\leq p\le q<\infty$
     and $B_2\geq  d(1- \frac{1}{p}- \frac{1}{q})$, set $$ D_2=d(1-\frac{1}{p}-\frac{1}{q}), \quad L_2=B_2-D_2.$$ 
    
    \begin{enumerate}
        \item If $B_1>  d(1- \frac{1}{p}- \frac{1}{q})$, set $D_1=\max\Big(0,d(1- \frac{1}{p}- \frac{1}{q})\big)$ and $ L_1=B_1-D_1$. 
        \item If $B_1\leq d(1- \frac{1}{p}- \frac{1}{q})$, set $D_1=B_1$ and $L_1=0$. 
      
    \end{enumerate}
  For the chosen values, $C(t)t^{L^T}=t^{\beta_{p,q}^T}$. Therefore,  setting $t= \norm{f |x|^A}_q^{\frac{1}{\beta_{p,q}^T +A}}$, we  obtain
            \begin{equation*}
            \qquad\qquad\norm{f}_q\lesssim \norm{f|x|^A}_q^{\frac{\beta^T_{p,q}}{A+\beta^T_{p,q}}} + \norm{f|x|^A}_q \norm{f|x|^A}_q^{-\frac{A}{A+\beta^T_{p,q}}} \approx  \norm{f|x|^A}_q^{\frac{\beta^T_{p,q}}{A+\beta^T_{p,q}}} .\end{equation*}
    \item Let $q<p$, $\max(q,p')\geq 2$, and $B_2 > d( 1 -\frac{1}{p}- \frac{1}{q})$.
    Set $$ \min(B_2,\frac{d}{p'})> D_2>d(1-\frac{1}{p}-\frac{1}{q}), \quad L_2=B_2-D_2.$$ 
    \begin{enumerate}
        \item If $B_1 \neq  d(1- \frac{1}{p}- \frac{1}{q})$, we proceed as in $(1)$ and obtain $C(t)t^{L^T}=t^{\beta_{p,q}^T}.$
        \item If $B_1= d(1- \frac{1}{p}- \frac{1}{q})$, we proceed as in $(1b)$ and obtain $$C(t)t^{L^T}=\begin{cases}
        t^{\beta_{p,q}^T},\quad\quad\quad \quad t<1, \\
        \log(t+1)^{\frac{1}{r}},\quad t\geq 1.
        \end{cases}$$
          The result follows by setting $t= \norm{f |x|^A}_q^{\frac{1}{\beta_{p,q}^T +A}}$.
        
\end{enumerate}
        \item Let $q<p$, $\max(q,p')<2$, and $B_2>d(1-\frac{1}{p}-\frac{1}{2})$.
        \begin{enumerate}
            \item If $B_1+d(\frac{1}{2}+ \frac{1}{p}-1)\neq 0$, we proceed as in $(2a)$  setting $q=2$ and obtain $C(t)t^{L^T}=t^{\beta_{p,2}^T} t^{d (\frac1q- \frac12)}.$ We derive  the result by setting
             $t= \norm{f |x|^A}_q^{\frac{1}{\beta_{p,2}^T+ d (\frac1q-\frac12) +A}}$.
            \item If $B_1+d(\frac{1}{2}+ \frac{1}{p}-1)=0$, we proceed as in $(2b)$ setting $q=2$ and obtain $$C(t)t^{L^T}=t^{d (\frac1q- \frac12)}\begin{cases}
        t^{\beta_{p,2}^T}, \qquad\qquad\quad t<1, \\
        \log(t+1)^{\frac{1}{2}-\frac1p}, \quad t\geq 1.
        \end{cases}$$ Here we let 
             $t= \norm{f |x|^A}_q^{\frac{1}{\beta_{p,2}^T+ d (\frac1q-\frac12) +A}}$ if $\norm{f |x|^A}_q<1$, and $$t=\norm{f |x|^A}_q^{\frac{1}{d (\frac1q-\frac12) +A}} \log\left(1+\norm{f |x|^A}_q  \right)^{-\frac{\frac 12- \frac1p}{d (\frac1q-\frac12) +A}}$$ if $\norm{f |x|^A}_q>1$.
      
        \end{enumerate}
        \item Let $q=\infty$. If $B_2>\frac{d}{p'}$, then we set $t=\norm{f |x|^A}_\infty^{\frac{1}{\beta_{p,\infty}^T +A}}$ in \eqref{eq:infinity}.
        \end{enumerate}

\end{proof}

\subsection{Proof of $\gtrsim$ part of Theorem \ref{theorem:GH}}
To show the optimality of 
$H$, we use the  functions 
of the type
$\displaystyle f(x)=\phi(x) \sum_{n \in \mathbb{Z}^d} c_n e^{2 \pi i \langle x, n\rangle}.$
 We now give some simple properties of such functions.
\begin{lemma}
\label{lemma:Y}
Let $1\leq p,q\leq \infty$ and $A,B    \geq 0$. Let $\widehat{\phi}$ be a smooth, non-negative function supported on the unit cube such that $\phi$ does not vanish on $1+ [0,1]^d$. Define
\begin{equation}
\label{eq:deff}
    f(x)=\phi(x) \sum_{n \in \mathbb{Z}^d} c_n e^{2 \pi i \langle x, n\rangle},
\end{equation}
where $c_n\in \mathbb{R}$ is only non-zero for finitely many $n$, and let $f_t(x)=f\left(\frac{x}{t}\right), t>0$.

Then
\begin{equation}
\label{eq:discret}
    \norm{f_t |x|^A} _q \approx t^{A+ \frac{d}{q}} \Big(\int_{[0,1]^d} \Big |\sum_{n\in \mathbb{Z}^d}  c_n e^{2 \pi i \langle n , x \rangle}  \Big|^q dx\Big)^{\frac{1}{q}}
\end{equation}
and
\begin{equation}
\label{eq:discret2}
     \norm{\widehat{f}_t |\xi|^B}_p \approx t^{\frac{d}{p'}} \Big(\sum_{n\in \mathbb{Z}^d} |c_n|^p\Big(\frac{1 +|n|}{t}\Big)^{pB}\Big)^{\frac{1}{p}}.
\end{equation}
\end{lemma}
\begin{proof}
Let us prove \eqref{eq:discret}.
Using periodicity and an appropriate change of variables, we deduce that
\begin{align*}
    \norm{f_t |x|^A}^q _q &= t^d \int_{\mathbb{R}^d} \Big ||tx|^A \phi(x) \sum_{n\in \mathbb{Z}^d} c_n e^{2 \pi i \langle x, n\rangle}\Big|^q dx\\
    & = t^d \int_{[0,1]^d}  \sum_{m\in \mathbb{Z}^d } \left ||tx+tm|^A \phi(x+m)\right|^q  \Big|\sum_{n\in \mathbb{Z}^d} c_n e^{2 \pi i \langle x,n\rangle}\Big|^q dx.
\end{align*}
To show
\begin{equation*}
   \sum_{m\in \mathbb{Z}^d} \left |tx+tm\right|^{Aq} |\phi(x+m)|^q \approx t^{Aq},\, x \in [0,1]^d,
\end{equation*}
we first note the lower bound follows from the fact that 
$\phi$ does not vanish on $1+ [0,1]^d$. Second, if $x>1$, the upper bound follows from the estimate $$\left |tx+tm\right|^A \lesssim 1+t^{A_2}(|x+m|)^{A_2}, $$ and the fact that $\phi$ decays faster than any power.
Third, if $t<1$, \begin{align*}
t^{A_1}&\gtrsim \Bigg(t^{A_2}\sum_{|t(x+m)|\geq 1} +  t^{A_1}\sum_{|t(x+m)|\leq 1}\Bigg) \left ||x+m|^{A_1+A_2} \phi(x+m)\right|^q\\
& \gtrsim \sum_{m\in \mathbb{Z}^d} \left ||tx+tm|^A \phi(x+m)\right|^q,
\end{align*} where 
 we again rely on the rapid decay of
$\phi$.

To prove \eqref{eq:discret2}, we observe that $$\widehat{f}_t(\xi)=t^d  \sum_{n\in \mathbb{Z}^d}  c_n \widehat{\phi}(t\xi-n).$$ Then, since
     $\widehat{\phi}$ is supported in the unit cube, the supports of $\widehat{\phi}(t\xi-n)$ are disjoint,  we conclude that 
     $$\norm{\widehat{f}_t |\xi|^B}^p_p=t^{dp}  \sum_{n\in \mathbb{Z}^d}  |c_n|^p  \norm{\widehat{\phi}(t\xi-n) |\xi|^B}_p^p \approx  t^{d(p-1)} \sum_{n\in \mathbb{Z}^d}  |c_n|^p \left(\frac{1 + |n|}{t}\right)^{pB}.$$

\end{proof}

To complete the proof of Theorem \ref{theorem:GH}, it remains only to choose convenient values for
 $(c_n)$ and $t$  for $f$ in \eqref{eq:deff} and $f_t$.
\begin{proof}[Proof of $\gtrsim$ part of Theorem \ref{theorem:GH}]\;\newline

{\it Case $x>1$.}
Let $M=t>1$. Let $\tilde{c}=(\tilde{c}_n)_{|n|\leq M}$ be such that $$C_{q,p}(M)\approx  \frac{ \left(\int_{[0,1]^d}\left | \sum_{|n|\leq M} \tilde{c}_n e^{2 \pi i \langle x, n\rangle}\right|^q\right)^{\frac{1}{q}}}{\left(\sum_{|n|\leq M} |\tilde{c}_n|^p(1+|n|)^{p\gamma}\right)^{\frac{1}{p}}},$$ where $C_{q,p}(M)$ is defined in Lemma \ref{lemma:X}. By Lemma \ref{lemma:Y}, there exists a multiple of $\tilde{c}$, which we call $c$, so that the function $f$ defined in \eqref{eq:deff} satisfies 
$$\norm{f_t}_q \approx M^{d(\frac{1}{q}+\frac{1}{p}-1)+B_1} C_{q,p}(M),\quad \norm{f_t |x|^A}_q \approx M^{d(\frac{1}{q}+\frac{1}{p}-1)+B_1+A_2} C_{q,p}(M),$$ together with the normalization condition  $$t^{-B_1}\left(\sum_{|n|\leq M} |c_n|^p (1+|n|)^{pB_1}\right)^{\frac{1}{p}}\approx \norm{\widehat{f_t} |\xi|^B}_p  =1.$$
\begin{enumerate}[label=(\Alph*)]
    \item Consider the following cases, corresponding to  items $(1),$ $(4a)$ and $(2a)$  of Theorem \ref{theorem:GH} in the regime $B_1<d(1-\frac 1p - \frac 1q)$: \begin{itemize}
        \item  $1\leq p\leq q < \infty$ and $B_1<d(1-\frac 1p - \frac 1q)$,
        \item $1\leq p\leq q = \infty$ and $B_1<d(1-\frac 1p)$,
        \item $1\leq q<p \leq \infty, \max(q,p')\geq 2$, and $B_1<d(1-\frac 1p - \frac 1q)$.
\end{itemize} We observe that, under these conditions, our goal is to prove the estimate (note that the first coordinate of $\beta_{p,q}$ is $0$) 
$$
H(x)\gtrsim x^{\frac{\beta^T_{p,q}}{A+{\beta^T_{p,q}}}}=1.
$$
    Here, by Lemma \ref{lemma:X}, $C_{q,p}(M)\approx M^{-B_1+ d(1-\frac1p-\frac 1q)}$. It follows that $$\norm{f_t}_q \approx 1 \mbox{ and } \norm{f_t |x|^A}_q \approx M^{A_2},$$ so that 
    $$H(x)\gtrsim 1,\qquad x>1.$$ In fact this inequality follows directly from the definition of $H$, which is an increasing function.
    \item Consider the following cases, corresponding to  item $(1)$ in the regime $B_1\geq d(1-\frac 1p - \frac 1q)$ and to items $(4a)$ and $(2a)$  in the regime $B_1>d(1-\frac 1p - \frac 1q)$:
    \begin{itemize}
        \item 
        $1\leq p\leq q < \infty$ and $B_1\geq d(1-\frac 1p - \frac 1q)$,
        \item  $1\leq p\leq q=\infty$ and $B_1> d(1-\frac 1p)$,
        \item  $1\leq q<p\leq \infty, \max(q,p')\geq 2$, and $B_1> d(1-\frac 1p - \frac 1q)$.
\end{itemize} Under these conditions, we want to show the lower bound $$H(x)\gtrsim x^{\frac{\beta^T_{p,q}}{A+{\beta^T_{p,q}}}}=x^{\frac{B_1+ d (\frac1p+\frac1q-1)}{A_2+B_1+ d (\frac1p+\frac1q-1)}}.$$
    In this case, $C_{q,p}(M)\approx 1$. Thus,
    $$\norm{f_t}_q \approx M^{d(\frac{1}{q}+\frac{1}{p}-1)+B_1} \mbox{ and } \norm{f_t |x|^A}_q \approx M^{d(\frac{1}{q}+\frac{1}{p}-1)+B_1+A_2},$$ which means that 
    $$H(x)\gtrsim x^{\frac{d(\frac{1}{q}+\frac{1}{p}-1)+B_1}{d(\frac{1}{q}+\frac{1}{p}-1)+B_1+A_2}},\qquad x>1.$$
    \item Consider the case $1\leq p \leq q=\infty$ and $B_1=d(1-\frac1p)$, corresponding to item $(4b)$. Note that $C_{q,p}(M)\approx \log(M+1)^{1-\frac1p}$. Then we have $$\norm{f_t}_q \approx\log(M+1)^{1-\frac1p} \mbox{ and } \norm{f_t |x|^A}_q \approx M^{A_2}\log(M+1)^{1-\frac1p},$$ which implies
    $$H(x)\gtrsim  \log(x+1)^{1-\frac{1}{p}},\qquad x>1.$$

     \item Consider the case  $1\leq q<p\leq\infty, \max(q,p')\geq 2$, and $B_1= d(1-\frac 1p - \frac 1q)$, corresponding to item (2)(b). Here $C_{q,p}(M)\approx \log(M+1)^{\frac1r}$. Thus $$\norm{f_t}_q \approx \log(M+1)^{\frac1r} \mbox{ and } \norm{f_t |x|^A}_q \approx M^{A_2}\log(M+1)^{\frac1r},$$ which means that 
    $$H(x)\gtrsim \log(x+1)^{\frac1r},\qquad  x>1.$$
     \item Consider the case $1\leq q<2<p\leq\infty$ and $B_1<d(1-\frac1p-\frac12)$, which corresponds to item $(3a)$ in the regime $B_1<d(1-\frac{1}{p}-\frac{1}{2})$. To obtain $$H(x)\gtrsim x^{\frac{\beta^T_{p,2}+d(\frac1q-\frac12)}{A+{\beta^T_{p,2}}+d(\frac1q-\frac12)}}=x^{\frac{d(\frac1q-\frac12)}{A+d(\frac1q-\frac12)}},$$
we note that,     in this case, $C_{q,p}(M)\approx M^{-B_1+ d(1-\frac1p-\frac12)}.$ It follows that $$\norm{f_t}_q \approx M^{d (\frac1q-\frac12)}\ \mbox{ and } \norm{f_t |x|^A}_q \approx M^{d (\frac1q-\frac12)+A_2},$$ and hence
    $$H(x)\gtrsim x^{\frac{d(\frac1q-\frac12)}{A_2+ d(\frac1q-\frac12)}}, \qquad x>1.$$
    \item Consider the case $1\leq q<2<p\leq\infty$ and $B_1=d(1-\frac1p-\frac12)$, corresponding to item $(3b)$. Note that $C_{q,p}(M)\approx \log(M+1)^{\frac12 - \frac1p}$. We have$$\norm{f_t}_q \approx M^{d (\frac1q-\frac12)}\log(M+1)^{\frac12- \frac1p}$$ and$$\norm{f_t |x|^A}_q \approx M^{d (\frac1q-\frac12)+A_2}\log(M+1)^{\frac12- \frac1p},$$ which yields 
    $$H(x)\gtrsim x^{\frac{d(\frac{1}{q}-\frac{1}{2})}{d(\frac{1}{q}-\frac{1}{2})+A_2}} \log(x+1)^{\frac{A_2 (\frac12-\frac1p)}{d(\frac{1}{q}-\frac{1}{2})+A_2}}, \qquad x>1.$$
    \item Consider the case $1\leq q<2<p\leq\infty$ and $B_1>d(1-\frac1p-\frac12)$, corresponding to $(3a)$ in the regime $B_1>d(1-\frac{1}{p}-\frac{1}{2})$. We would like to verify that $$H(x)\gtrsim x^{\frac{\beta^T_{p,2}+d(\frac1q-\frac12)}{A+{\beta^T_{p,2}}+d(\frac1q-\frac12)}}=x^{\frac{d(\frac{1}{q}+\frac{1}{p}-1)+B_1}{d(\frac{1}{q}+\frac{1}{p}-1)+B_1+A_2}}.$$
    Indeed, here $C_{q,p}(M)\approx 1$, which means that $$\norm{f_t}_q \approx M^{d (\frac1q+\frac1p-1)+B_1}\ \mbox{ and } \norm{f_t |x|^A}_q \approx M^{d (\frac1q+\frac1p-1)+B_1+A_2},$$ $$H(x)\gtrsim x^{\frac{d(\frac{1}{q}+\frac{1}{p}-1)+B_1}{d(\frac{1}{q}+\frac{1}{p}-1)+B_1+A_2}},\qquad x>1.$$ 
\end{enumerate}

{\it Case $x<1$. } Similarly as in the case $x>1$, for every $M\geq 1>t$, there exists a sequence $c_n$  such that 
$$\norm{f_t}_q \approx t^{d(\frac{1}{q}+\frac{1}{p}-1)+B_2} C_{q,p}(M), \quad \norm{f_t |x|^A}_q \approx t^{d(\frac{1}{q}+\frac{1}{p}-1)+B_2+A_1} C_{q,p}(M),$$ and the normalization condition $\norm{\widehat{f_t} |\xi|^B}_p  =1$ is satisfied.
We observe  that here  our goal  is to show that in  cases $(1)$--$(4)$ of Theorem \ref{theorem:GH}
\begin{equation}
\label{eq:Hfrombelow}    
H(x)\gtrsim x^{\frac{B_2 + d(\frac1q+\frac1p-1)}{A_1+B_2 + d(\frac1q+\frac1p-1)}},\end{equation} and, in those cases not covered by items $(1)$--$(4)$, that $H(x)=\infty$. 
Assume first that \begin{equation}
\label{eq:assumption}
    d(\frac{1}{q}+\frac{1}{p}-1)+B_2+A_1>0.
\end{equation}
Let $$t=\left(\frac{x}{C_{q,p}(M)}\right)^{\frac{1}{{d(\frac{1}{q}+\frac{1}{p}-1)+B_2+A_1}}}\lesssim 1.$$ Then we have
\begin{equation}
\label{vsp}    
H(x)\gtrsim x^{\frac{B_2 + d(\frac1q+\frac1p-1)}{A_1+B_2 + d(\frac1q+\frac1p-1)}}C_{q,p}(M)^{\frac{A_1 }{A_1+B_2 + d(\frac1q+\frac1p-1)}}, \; x\leq 1\leq M.\end{equation}
Hence, taking $M=1$, we deduce \eqref{eq:Hfrombelow} for 
 $x\leq 1.$

In those cases not mentioned in items (1)--(4) of Theorem \ref{theorem:GH} we can obtain a stronger inequality than \eqref{eq:Hfrombelow} by letting $M\to \infty$. Indeed, one of the following must hold: \begin{itemize}
\item $1\leq p,q\leq \infty$ and $B_2< d(1-\frac1p-\frac1q)$,
    \item $1\leq q<p \leq \infty, \max(q,p')\geq 2$, and $B_2=  d(1-\frac 1p - \frac 1q)$, \item $1\leq q<2<p \leq \infty $ and $B_2\leq  d(1-\frac 1p - \frac 12)$,
     \item $1<p\le q=\infty$
    and $B_2= d(1-\frac1p)$.
\end{itemize} Then, by Lemma \ref{lemma:X}, we know that $C_{q,p}(M)\to \infty$ as $M\to \infty$, so inequality \eqref{vsp} implies that $$H(x)=\infty.$$

Finally, if \eqref{eq:assumption} fails, that is, if $d(\frac{1}{q}+\frac{1}{p}-1)+B_2+A_1\leq 0$, there exists $\tilde{B_2}$
satisfying
$B_2 <\tilde{B_2}<d(1-\frac{1}{q}-\frac{1}{p})$ and
$d(\frac{1}{q}+\frac{1}{p}-1)+\tilde{B_2}+A_1> 0$, which is  assumption \eqref{eq:assumption} with $\tilde{B_2}$. Hence, the preceding results imply that $$H_{A,(B_1,\tilde{B_2})}=\infty.$$ Since $H$ is decreasing with respect to $B_2$, we conclude that $H_{A,(B_1,{B_2})}=\infty$.
\end{proof}

\begin{proof}[Proof of
Corollary \ref{corollary-broken}]

  Inequality \eqref{vspom1}  is equivalent to $$\norm{f}_q \lesssim \norm{\widehat{f}|\xi|^B}_{q'} \left(\frac{\norm{f |x|^A}_q}{\norm{\widehat{f} |\xi|^B}_
{q'}}\right)^{\frac{B^T}{A+B^T}},$$ 
which follows from Theorem \ref{theorem:GH}.
\end{proof}

\section{Symmetric uncertainty inequalities
}
\label{section:broken}
In this section, we prove Theorem \ref{theorem:broken}, which characterizes the parameters $A,B,p,p'$ for which inequality \eqref{ineq:uncer--}, that is, \begin{equation}   \norm{f}_{p}\norm{\widehat{f}}_{p'} \lesssim \norm{|x|^A f}_p \norm{|\xi|^B \widehat{f}}_{p'}
    \label{ineq:uncer--.}
\end{equation} holds.
\subsection{Proof of sufficiency  in Theorem \ref{theorem:broken}}

\begin{proof}
 Before we begin, we note that the sufficiency of the third condition in Theorem \ref{theorem:broken} 
  follows from the second condition as follows. Assume that $d(\frac{1}{2}-\frac{1}{p}) =A_1 <B_2$ and set $A':=(B_2,A_2)$. Then, assuming that condition (2) is sufficient, we have
 $$\norm{f}_{p} \norm{\widehat{f}}_{p'} \lesssim \norm{f|x|^{A'}}_p \norm{\widehat{f} |x|^B}_{p'} \leq \norm{f|x|^{A}}_p \norm{\widehat{f} |x|^B}_{p'}.$$

 We may now assume that $A_1 \neq d(\frac{1}{2}-\frac{1}{p})$. The proof of the remaining items follows directly from Theorem \ref{theorem:GH}. Indeed,  applying Theorem \ref{theorem:GH}$(1a)$ (or, if $p=\infty$, $(4a)$)  to the function $\displaystyle\frac{f
}{\norm{\widehat{f} |\xi|^B}_{p'}}$,  
we have
  $$\frac{\norm{f}_p}{
\norm  {\widehat{f} |\xi|^B}_{p'}}
  \lesssim  H_1\Bigg(\frac{\norm{f |x|^A}_p}{\norm{\widehat{f} |\xi|^B}_{p'}}\Bigg),$$
where
$$H_1(x)= x^{\frac{B^T}{A+ B^T}}.$$  Similarly, $$\frac{\norm{\widehat{f} }_{p'}}{\norm{{f}|x|^A}_{p}} \lesssim H_2\Bigg(\frac{\norm{\widehat{f} |\xi|^B}_{p'}}{\norm{f |x|^A}_p}\Bigg),$$ for a suitable $H_2$.  Therefore, inequality \eqref{ineq:uncer--} holds provided that $$H_1(x) H_2(x^{-1}) \lesssim 1, \mbox{ for } x\in \mathbb{R}_+.$$

{\it Case 1.}
Let $p=2$. Then, by Theorem \ref{theorem:GH} $(1a)$, $H_2(x)= x^{\frac{A^T}{B+ A^T}}$, and
 $$H_1(x)H_2(x^{-1})= x^{\frac{B^T}{A+B^T}} (\frac{1}{x})^{\frac{A^T}{B+A^T}} =
  x^{\frac{B^T-A}{A+B^T}}\le1. $$
  
  {\it Case 2.}
Let $p>2$. Recalling that 
$A_1 + d(\frac{1}{p} -\frac{1}{2}) \neq 0,$
and setting
    $$\alpha_{p,2}:= \left(\max\Big(A_1 + d(\frac{1}{p} -\frac{1}{2}),0\Big), A_2 +d (\frac{1}{p} -\frac{1}{2})\right),$$
Theorem \ref{theorem:GH} $(3a)$, applied to  $\displaystyle\frac{\widehat{f}
}{\norm{f |x|^A}_p}$,
yields 

  $$H_2(x)=x^{\frac{\alpha_{p,2}^T+
d(\frac{1}{p'} -\frac{1}{2})}{B+\alpha_{p,2}^T +d(\frac{1}{p'} -\frac{1}{2})}}.$$

  {\it Case 2.1.} Let $A_1>  d(\frac{1}{2}-\frac{1}{p})$. In this case
    $\alpha_{p,2}= A+ d(\frac{1}{p} -\frac{1}{2})$
    and
$$H_1(x)H_2(x^{-1}) = x^{\frac{B^T-A}{A+B^T}} \leq 1,$$
because   
 $B_2 \geq A_1$ and $ A_2 \geq B_1$.

  {\it Case 2.2.} Let $A_1<  d(\frac{1}{2}-\frac{1}{p})$. In this case
    $\alpha_{p,2}= (0, A_2+ d(\frac{1}{p} -\frac{1}{2}))
    $
    and
$H_1(x)H_2(x^{-1}) = x^{C}$,
where $C=(C_1,C_2)$ with  $$C_1=\frac{B_2}{A_1+B_2}
-
\frac{
d(\frac{1}{p'} -\frac{1}{2})}{B_2+d(\frac{1}{p'} -\frac{1}{2})}
$$
and  $$C_2=\frac{B_1}{A_2+B_1}
-
\frac{
A_2}{B_1+A_2}
.$$ We conclude the proof by noting that  $x^C\leq 1$, because both
$C_1\geq 0$ (due to
$A_1 d(\frac{1}{p'} -\frac{1}{2})\leq B_2^2$)
and 
$C_2\leq0$ (since
$A_2 \geq B_1)$ hold.
\end{proof}

\subsection{Proof of necessity  in Theorem \ref{theorem:broken}}

For the proof of necessity, we first establish two simple preliminary results. The first serves as a substitute for the usual homogeneity relation $\displaystyle\norm{f( x/\lambda) |x|^{\alpha}}_p = \lambda^{\alpha +\frac{d}{p}} \norm{f( x) |x|^{\alpha}}_p$ in the case of broken power weights. The second provides some properties of the functions defined in \eqref{eq:counterex}, which will be our source of counterexamples to \eqref{ineq:uncer--}.

 \begin{lemma}
 \label{lemma:escaling}
 Assume that inequality \eqref{ineq:uncer--} holds for some $A$ and $B$. Then,
\begin{itemize}
\item For any $\lambda>0$, we have
\begin{equation}
\begin{aligned} 
     &\left( \lambda ^{-pA_1} \int_{|x|\leq \lambda} |f(x)|^p |x|^{pA_1} dx + \lambda ^{-pA_2} \int_{|x|\geq \lambda} |f(x)|^p |x|^{pA_2} dx\right)^{\frac{1}{p}}  \\
   &\times\left( \lambda ^{p'B_1} \int_{|\xi|\leq \lambda^{-1}} |\widehat{f}(\xi)|^{p'} |\xi|^{p'B_1} d\xi + \lambda ^{p'B_2} \int_{|\xi|\geq \lambda^{-1}} |\widehat{f}(\xi)|^{p'} |\xi|^{p'B_2} d\xi\right)^{\frac{1}{p'}} \\
&\qquad\qquad\qquad\qquad\qquad\qquad\qquad\qquad\qquad\qquad\qquad\lesssim \norm{f}_p\norm{\widehat{f}}_{p'};\label{eq:lambdaeq}
   \end{aligned}
    \end{equation}
    \item $A_1\leq B_2$ and $B_1\leq A_2;$

    \item If $A_2=B_1$, then inequality \eqref{ineq:uncer--} also holds for $A=(A_2,A_2)$ and $B=(B_1,B_1);$
    \item If $A_1=B_2$, then inequality \eqref{ineq:uncer--} also holds for $A=(A_1,A_1)$ and $B=(B_2,B_2)$.
\end{itemize}
 \end{lemma}
\begin{proof}

For a function $f$ and $\lambda>0$, define
$$f_\lambda(x)=\lambda^d f(\lambda x).$$
Then, $\widehat{f_\lambda}(\xi)=\widehat{f}(\xi/\lambda)$. Thus,
$\norm{f_\lambda}_p = \lambda ^{d/p'} \norm{f}_p$ and 
$\norm{\widehat{f_\lambda}}_{p'} = \lambda ^{d/p'} \norm{\widehat{f}}_{p'}.$
Also, with the usual modification for $p=\infty,$
$$\lambda^{-\frac{dp}{p'}}\norm{f_\lambda x^A}_p^p =   
\int_{|x|\leq \lambda} |f(x)|^p \Big(\frac{|x|}{\lambda}\Big)^{pA_1} dx + 
\int_{|x|\geq \lambda} |f(x)|^p \Big(\frac{|x|}{\lambda}\Big)^{pA_2}dx,$$
$$\lambda^{-d}\norm{\widehat{f_\lambda} \xi ^B}_{p'}^{p'} =   
\int_{|\xi|\leq \lambda^{-1}} |\widehat{f}(\xi)|^{p'} (\lambda|\xi|)^{p'B_1} d\xi + 
\int_{|\xi|\geq \lambda^{-1}} |\widehat{f}(\xi)|^{p'} (\lambda|\xi|)^{p'B_2}d\xi,$$ 
 whence \eqref{eq:lambdaeq} follows.

To prove the second item, we note that if $f$ is a Schwartz function 

\begin{equation}
    \label{eq:vsp1}\norm{f}_p \frac{\norm{f_\lambda x^A}_p}{\norm{f_\lambda}_p} \approx \lambda^{-A_1} \left( \int_{|x|\leq \lambda} |f(x)|^p |x|^{pA_1} dx\right)^{\frac{1}{p}} \mbox{ for } \lambda \mbox{ large enough},
\end{equation}
because
$\left(\int_{|x|\geq \lambda} |f(x)|^p |x|^{pA_2} dx\right)^{\frac{1}{p}}$ decays faster than any polynomial. Besides, by observing that $$\lim_ {\lambda \to \infty} \lambda ^{B_1} \left(\int_{|\xi|\leq \lambda^{-1}} |\widehat{f}(\xi)|^{p'} |\xi|^{p'B_1} d\xi\right)^{\frac{1}{p'}} <\infty,$$ we conclude that 
\begin{equation}
    \label{eq:vsp2}
    \norm{\widehat{f}}_{p'}\frac{\norm{\widehat{f_\lambda} x^B}_p}{\norm{\widehat{f_\lambda}}_{p'}} \approx \lambda^{B_2} \left( \int_{|x|\geq \lambda^{-1}} |\widehat{f}(\xi)|^{p'} |\xi|^{p'B_2} dx\right)^{\frac{1}{p'}},
\end{equation} for $\lambda$ large enough. The restriction $B_2\geq A_1$ follows by multiplying \eqref{eq:vsp1} and \eqref{eq:vsp2} and letting $\lambda \to \infty$. An analogous reasoning shows that
  $A_2\geq B_1$ must hold.

Finally, note that if $B_2=A_1$, multiplying \eqref{eq:vsp1} and \eqref{eq:vsp2} we deduce that
$$ \norm{f}_p \norm{\widehat{f}}_{p'}\lesssim \left( \int_{|x|\leq \lambda} |f(x)|^p |x|^{pA_1} dx\right)^{\frac{1}{p}}\left( \int_{|x|\geq \lambda^{-1}} |\widehat{f}(\xi)|^{p'} |\xi|^{p'B_2} dx\right)^{\frac{1}{p'}},$$ whence the forth item follows by letting $\lambda \to \infty$. The third item is proved in the same manner.
\end{proof}
\begin{lemma}
\label{lemma:auxi}Let $2\leq p \leq \infty$ and $A,B\geq 0$. 
     Let  $0\leq \phi \leq 1$ be  a smooth function supported on $|x|\leq \frac{1}{2}$ such that $\widehat{\phi}$ does not vanish on the cube $1+ [0,1]^d$. Define
     \begin{equation}
     \label{eq:counterex}
         f(x)=f_1(x) +f_2(x):= N^d \phi(Nx) + \sum_{2 \leq |n|\leq M} \varepsilon _n c_n \phi(x-n),
     \end{equation} where  $\lambda\in \mathbb{R},$ $M\in \mathbb{N}$, $c_n\in \mathbb{R}$, and $\varepsilon_n$ is either $1$ or $-1$.
     Then,
   \begin{align}
       \label{eq:lemma4.3-1}
       \norm{f}_p \approx N^{\frac{d}{p'}}+ \Big(\sum_{2\leq |n|\leq M} |c_n|^p \Big)^{\frac{1}{p}}
   \end{align}  
 and  \begin{align} 
     \norm{ f \big|\frac{x}\lambda\big|^A}_p 
     &=\Bigg( \lambda ^{-pA_1} \int_{|x|\leq \lambda} |f(x)|^p |x|^{pA_1} dx + \lambda ^{-pA_2} \int_{|x|\geq \lambda} |f(x)|^p |x|^{pA_2} dx\Bigg)^{\frac{1}{p}}
     \nonumber \\
    \label{eq:lemma4.3-2}     &\approx \lambda^{-A_1} \Bigg( N^{\frac{d}{p'}-A_1} \norm{ |x|^{A_1} \phi \mathbbm{1}_{|x|\leq N\lambda}}_p + \Big(\sum_{2\leq |n|\leq \lambda} |c_n|^p |n|^{pA_1}\Big)^{\frac{1}{p}}\Bigg)
         \\
         &+\lambda^{-A_2} \Bigg( N^{\frac{d}{p'}-A_2} \norm{ |x|^{A_2} \phi \mathbbm{1}_{|x|\geq N\lambda}}_p + \Big(\sum_{\lambda < |n|\leq M} |c_n|^p |n|^{pA_2}\Big)^{\frac{1}{p}} \Bigg).
     \end{align}

Moreover, for 
$\mathbbm{E}_{\varepsilon} f(\varepsilon_1, \dots, \varepsilon_M)=\frac{1}{2^M} \sum_{\varepsilon_1=\pm 1, \dots, \varepsilon_M=\pm 1} f(\varepsilon_1, \dots, \varepsilon_M)$, we have 
   \begin{align}
\label{eq:lemma4.3-3}  N ^{d}+\mathbbm{E}_{\varepsilon} \norm{\widehat{f}}^{p'}_{p'} \gtrsim \Big(\sum_{2 \leq |n|\leq M} |c_n|^2 \Big)^{\frac{p'}{2}}
     \end{align}
     and
     \begin{align}
\label{eq:lemma4.3-4}    
&\left(\mathbbm{E}_{\varepsilon} \norm{\widehat{f} |\lambda \xi|^{B}}^{p'}_{p'}\right)^{\frac{1}{p'}}\\
     &=\left( \lambda ^{p'B_1} \int_{|\xi|\leq \lambda^{-1}} \mathbbm{E}_{\varepsilon}|\widehat{f}(\xi)|^{p'} |\xi|^{p'B_1} d\xi + \lambda ^{p'B_2} \int_{|\xi|\geq \lambda^{-1}} \mathbbm{E}_{\varepsilon}|\widehat{f}(\xi)|^{p'} |\xi|^{p'B_2} d\xi\right)^{\frac{1}{p'}}\\
      &\lesssim \lambda^{B_1}\Bigg(N^{\frac{d}{p'}+B_1} \norm{|\xi|^{B_1}\widehat{\phi} \mathbbm{1}_{|\xi|\leq (N\lambda)^{-1}}}_{p'} +  \Big(\sum_{2 \leq |n|\leq M} |c_n|^2 \Big)^{\frac12}\norm{|\xi|^{B_1}\widehat{\phi} \mathbbm{1}_{|\xi|\leq \lambda^{-1}}}_{p'}\Bigg)\\
      &+\lambda^{B_2} \Bigg(N^{\frac{d}{p'}+B_2} \norm{|\xi|^{B_2}\widehat{\phi} \mathbbm{1}_{|\xi|\geq (N\lambda)^{-1}}}_{p'} + \Big(\sum_{2 \leq |n|\leq M} |c_n|^2 \Big)^{\frac12}\norm{|\xi|^{B_2}\widehat{\phi} \mathbbm{1}_{|\xi|\geq \lambda^{-1}}}_{p'}\Bigg).
  \end{align}
\end{lemma}
\begin{proof}
    The proof of \eqref{eq:lemma4.3-1}  follows by noting that $|f|^p =|f_1|^p + |f_2|^p$ and applying Lemma \ref{lemma:Y} to $f_2$. 
    Similarly, for \eqref{eq:lemma4.3-2} we have \begin{align*}
     \norm{f_1 |x/\lambda|^A}_p&\approx
         \lambda^{-A_1}  N^{\frac{d}{p'}-A_1} \norm{ |x|^{A_1} \phi \mathbbm{1}_{|x|\leq N\lambda}}_p\\ &+\lambda^{-A_2}  N^{\frac{d}{p'}-A_2} \norm{ |x|^{A_2} \phi \mathbbm{1}_{|x|\geq N\lambda}}_p,
    \end{align*}
    and, as in Lemma \ref{lemma:Y},
    \begin{align*}
        \norm{f_2 |x/\lambda|^A}_p&\approx \lambda^{-A_1}  \Big(\sum_{2\leq |n|\leq \lambda} |c_n|^p |n|^{pA_1}\Big)^{\frac{1}{p}}+ \lambda^{-A_2}\Big(\sum_{\lambda < |n|\leq M} |c_n|^p |n|^{pA_2}\Big)^{\frac{1}{p}}.
    \end{align*}
   To obtain \eqref{eq:lemma4.3-3}, we  have, once again by Lemma \ref{lemma:Y},
    $$\norm{\widehat{f}}^{p'}_{p'}+  \norm{\widehat{f_1}}_{p'}^{p'}\gtrsim  \Big\|\sum_{2\leq |n|\leq M} \varepsilon_n c_n e^{2 \pi i \langle n,x\rangle}
\Big\|^{p'}_{p'} . $$
    The result follows by observing that  the Khinchin inequality implies
    $$\mathbbm{E}_{\varepsilon}
\Big\|\sum_{2\leq |n|\leq M} \varepsilon_n c_n e^{2 \pi i \langle n,x\rangle}
\Big\|^{p'}_{p'}
  \approx \Big(\sum_{2\leq |n|\leq M}  |c_n|^2\Big)^{\frac{p'}{2}},  $$
    and that  $\norm{\widehat{f_1}}_{p'} \approx N^{\frac{d}{p'}}$.

    For \eqref{eq:lemma4.3-4}, similarly to \eqref{eq:lemma4.3-2}, we have  \begin{align*}\norm{\widehat{f_1} |\xi\lambda|^{B}}_{p'}&\approx\lambda^{B_1}N^{\frac{d}{p'}+B_1} \norm{|\xi|^{B_1}\widehat{\phi} \mathbbm{1}_{|\xi|\leq (N\lambda)^{-1}}}_{p'}\\
    &+\lambda^{B_2} N^{\frac{d}{p'}+B_2} \norm{|\xi|^{B_2}\widehat{\phi} \mathbbm{1}_{|\xi|\geq (N\lambda)^{-1}}}_{p'}. \end{align*}
For $\widehat{f_2},$ the Khinchin inequality gives
    $$\mathbbm{E}\norm{\widehat{f_2} |\xi\lambda|^{B}}_{p'}^{p'} \approx \Big(\sum_{2 \leq |n|\leq M} |c_n|^2\Big)^{\frac{p'}{2}} \norm{\widehat{\phi} |\xi\lambda|^B}_{p'}^{p'}, $$ whence the result follows.
\end{proof}
We are now in a position to prove the necessity part in Theorem \ref{theorem:broken}. The proof proceeds by selecting suitable values of $N,\lambda,$ and $(c_n)$.
        \begin{proof}[Proof of necessity in Theorem \ref{theorem:broken}]

   {\it Step 1.} To begin with, the conditions $B_2\geq A_1$ and $A_2\geq B_1$ follow from the second item of Lemma \ref{lemma:escaling}.

  {\it Step 2.} We now show that if $A_2 \leq d(\frac{1}{2}- \frac{1}{p})$, then \eqref{ineq:uncer--} cannot hold.
     Let $M,N \in \mathbb{N}$ .
     Since $A_2 \leq d(\frac{1}{2}- \frac{1}{p})$, there exists a sequence $c_n$ such that $\sum_{1 \leq |n| < \infty} |c_n| ^p |n| ^{pA_2}=:S^p < \infty$ and $\sum_{1 \leq |n| < \infty } |c_n|^2 = \infty$, for example, $c_n=\frac{1}{\sqrt{|n|^d \log (1+|n|)}}$.
     
    From inequality \eqref{ineq:uncer--} we deduce 
    $$\mathbbm{E}\left(\norm{f}^{p'}_p \norm{\widehat{f}}^{p'}_{p'}\right) \lesssim \mathbbm{E} \left(\norm{f |x|^A}^{p'}_p \norm{\widehat{f} |\xi |^B}^{p'}_{p'}\right).$$ Observing that $\norm{f}_p$ and $\norm{f |x|^A}_p$ do not depend on $\varepsilon$, we obtain 
    $$\norm{f}^{p'}_p\mathbbm{E} \norm{\widehat{f}}^{p'}_{p'} \lesssim \norm{f |x|^A}^{p'}_p\mathbbm{E}  \norm{\widehat{f} |\xi |^B}^{p'}_{p'}.$$ Dividing this last expression
    by $\left(\sum_{2\leq |n|\leq M} |c_n|^2\right)^{\frac12}$ and using
\eqref{eq:lemma4.3-3} 
and 
\eqref{eq:lemma4.3-4} from 
Lemma \ref{lemma:auxi}, we conclude, as $M\to \infty$, that 
$$\norm{f}_p \lesssim \norm{f|x|^A}_p \norm{|\xi|^{B}\widehat{\phi} }_{p'}.$$
In light of  \eqref{eq:lemma4.3-1} and \eqref{eq:lemma4.3-2}, 
we see
 that 
     \begin{align*}
         N^{\frac{d}{p'}} \lesssim \left(1+N^{\frac{d}{p'}-A_1} \right)\norm{|\xi|^{B}\widehat{\phi}}_{p'}.
     \end{align*}  Letting $N\to \infty$, we see that 
   this is only possible if $A_1=0$. 
    
     {\it Step 3.} If $A_2> d(\frac{1}{2}- \frac{1}{p})$ and $A_1< d(\frac{1}{2}- \frac{1}{p})$, we show that $B_2^2 \geq A_1d( \frac{1}{2}- \frac{1}{p})$. 
For $N$ large enough, define $\lambda = N ^{\frac{A_1+ B_2}{d(\frac{1}{2}- \frac{1}{p}) -A_1}}$ and $M=\lambda$.
Choose $c_n$ such that the Hölder inequality $$\Big(\sum_{2 \leq |n| \leq \lambda} |c_n|^2 \Big)^{\frac{1}{2}} \le \Big(\sum_{2 \leq |n| \leq \lambda} |c_n|^p |n| ^{p A_1}\Big)^{\frac{1}{p}} \Big(\sum_{2 \leq |n| \leq \lambda} |n|^{-\frac{A_1}{\frac{1}{2}- \frac{1}{p}}}\Big)^{\frac{1}{2}-\frac{1}{p}}$$ becomes an equality and normalize it so that
 $$\Big(\sum_{2 \leq |n| \leq \lambda} |c_n|^p |n| ^{p A_1}\Big)^{\frac{1}{p}} = N^{\frac{d}{p'}-A_1}.$$
  This implies that
 $$\Big(\sum_{2 \leq |n| \leq \lambda} |c_n|^2 \Big)^{\frac{1}{2}} \approx N ^{\frac{d}{p'}-A_1} \lambda ^{-A_1 + d(\frac{1}{2} - \frac{1}{p})}=N^{B_2+ \frac{d}{p'}}.$$
Besides, Lemma \ref{lemma:auxi} gives,  for $N$ large enough, that 
$$\norm{f}_p \gtrsim N^{\frac{d}{p'}};$$
  \begin{align*}
         &\Big( \lambda ^{-pA_1} \int_{|x|\leq \lambda} |f(x)|^p |x|^{pA_1} dx + \lambda ^{-pA_2} \int_{|x|\geq \lambda} |f(x)|^p |x|^{pA_2} dx\Big)^{\frac{1}{p}}\\
         &\approx
         \lambda^{-A_1} \Big( N^{\frac{d}{p'}-A_1} \norm{ |x|^{A_1} \phi \mathbbm{1}_{|x|\leq N\lambda}}_p + \Big(\sum_{2\leq |n|\leq \lambda} |c_n|^p |n|^{pA_1}\Big)^{\frac{1}{p}}\Big)
         \\
         &+\lambda^{-A_2} \Big( N^{\frac{d}{p'}-A_2} \norm{ |x|^{A_2} \phi \mathbbm{1}_{|x|\geq N\lambda}}_p + \Big(\sum_{\lambda < |n|\leq M} |c_n|^p |n|^{pA_2}\Big)^{\frac{1}{p}} \Big)\\
         & \approx  \lambda^{-A_1} N^{\frac{d}{p'}-A_1}
     \end{align*}
    (here we have used that $\phi$ decays faster than any power);
  $$\Big(\mathbbm{E}\norm{\widehat{f}}^{p'}_{p'}\Big)^{\frac{1}{p'}} \gtrsim  N ^{B_2+\frac{d}{p'}};$$
  and, setting $\|c\|_{\ell_2}:=\Big(\sum_{2 \leq |n|\leq M} |c_n|^2 \Big)^{\frac12}
$,
  \begin{align*}
      &\Big( 
\int_{|\xi|\leq \lambda^{-1}} \mathbbm{E}|\widehat{f}(\xi)|^{p'} (\lambda|\xi|)^{p'B_1} d\xi +  \int_{|\xi|\geq \lambda^{-1}} \mathbbm{E}|\widehat{f}(\xi)|^{p'} (\lambda|\xi|)^{p'B_2} d\xi\Big)^{\frac{1}{p'}}\quad\\
      &\lesssim \lambda^{B_1}\Big(N^{\frac{d}{p'}+B_1} \norm{|\xi|^{B_1}\widehat{\phi} \mathbbm{1}_{|\xi|\leq (N\lambda)^{-1}}}_{p'} +  \|c\|_{\ell_2}\norm{|\xi|^{B_1}\widehat{\phi} \mathbbm{1}_{|\xi|\leq \lambda^{-1}}}_{p'}\Big)\\
      &+\lambda^{B_2} \Big(N^{\frac{d}{p'}+B_2} \norm{|\xi|^{B_2}\widehat{\phi} \mathbbm{1}_{|\xi|\geq (N\lambda)^{-1}}}_{p'} + 
\|c\|_{\ell_2}\norm{|\xi|^{B_2}\widehat{\phi} \mathbbm{1}_{|\xi|\geq \lambda^{-1}}}_{p'}\Big)\\
      &\lesssim \lambda^{B_1}\left(N^{\frac{d}{p'}+B_1} (N\lambda)^{-B_1 - \frac{d}{p'}} +  N^{B_2+\frac{d}{p'}}\lambda^{-B_1-\frac{d}{p'}}\right)+\lambda^{B_2} N^{\frac{d}{p'}+B_2} \\
      &\approx N^{B_2 + \frac{d}{p'}} \lambda^{B_2}.
  \end{align*}

     Hence, inequality \eqref{eq:lambdaeq} implies
     $$N^{\frac{d}{p'}}N^{B_2+\frac{d}{p'}} \lesssim  N^{B_2 + \frac{d}{p'}} \lambda^{B_2}  \lambda^{-A_1} N^{\frac{d}{p'}-A_1},$$ that is,
   $$1 \lesssim \lambda^{B_2-A_1}N ^{-A_1}= N ^{\frac{(B_2-A_1) (A_1+ B_2)}{d(\frac{1}{2}- \frac{1}{p}) -A_1}-A_1}.$$ Letting $N \to \infty$ we conclude that
   $$0 \leq  \frac{(B_2-A_1) (A_1+ B_2)}{d(\frac{1}{2}- \frac{1}{p}) -A_1}-A_1.$$
   Equivalently,
   $B_2^2 -A_1^2 \geq -A_1^2 +A_1d( \frac{1}{2}- \frac{1}{p}),$ whence the result follows.
   
  {\it Step 4.}
    If $B_2=A_1=d(\frac{1}{2}- \frac{1}{p})$,
    the fourth item of Lemma \ref{lemma:escaling} shows that, should \eqref{ineq:uncer--} hold, it would also have to hold for $A=B=d(\frac{1}{2}- \frac{1}{p})$, which is impossible by Step 1.
     \end{proof}

  

\section{Heisenberg type uncertainty for special
Hadamard operators}\label{section: special}
As discussed in the introduction, there has been growing interest in establishing uncertainty principles for general operators that possess properties analogous to those of the Fourier transform. For example, in \cite{dias} the authors obtain several Heisenberg-type
uncertainty inequalities for operators $T$ satisfying \begin{equation}
\label{eq:jointboundedness}
    \norm{Tf}_\infty\leq \norm{f}_1 \quad \mbox{and}\quad\norm{Tf}_2 \leq \norm{f}_2
\end{equation} as well as the property $T^*=T^{-1}$.
We will refer to such operators as \textit{special Hadamard operators}. 

In what follows, 
$H(\cdot)$ denotes the function introduced in Theorem \ref{theorem:GH}.
For simplicity, in this section we consider only classical (non-broken) weights $|x|^\alpha$ and 
$|\xi|^\beta$, that is, the case
$\alpha=A_1=A_2$
and 
$\beta=B_1=B_2$.

The main goal of this section is to show that {\it some, but not all,} of the upper bounds of $H$ in Theorem \ref{theorem:GH} remain valid for all special Hadamard operators. In particular,  
the estimate 
\begin{equation}
\label{eq:uncerxiaoprim-}
    \norm{Tf}_q\lesssim \norm{Tf |x| ^{\alpha}}^\frac{\beta+d(\frac1p+\frac1q-1)}{\alpha+\beta+d(\frac1p+\frac1q-1)}_q \norm{{f} |\xi| ^\beta}^\frac{\alpha}{\alpha+\beta+d(\frac1p+\frac1q-1)}_{p}
\end{equation} holds for all special
Hadamard operators only for certain parameter values within the range specified in Theorem \ref{theo:mikxiao}.

First, we note that for small $\beta$ an analogue of Theorem \ref{theorem:GH} is valid for any operator satisfying \eqref{eq:jointboundedness}. 

\begin{proposition}
\label{prop:general}
    Let $1\leq p\leq \infty$, $1\leq q< \infty$, $\alpha>0$, and $0<\beta<d/p'$.
    Then
$$\sup _T \sup _{ \norm{Tf |x|^\alpha}_q\leq x, \norm{f |\xi|^\beta}_p\leq 1}  {\norm{Tf}_q}\approx H(x),$$ where the supremum is taken over all operators satisfying \eqref{eq:jointboundedness}.

In particular,
 under the conditions on $\beta, p,q$   stated in Theorem
    \ref{theo:mikxiao}, inequality 
    \eqref{eq:uncerxiaoprim-}
    holds for any operator satisfying
    \eqref{eq:jointboundedness}.

\end{proposition}
\begin{proof}Clearly, it suffices to show that
 $$\sup _T \sup _{ \norm{Tf |x|^\alpha}_q\leq x, \norm{f |\xi|^\beta}_p\leq 1}  {\norm{Tf}_q} \lesssim H(x)=x^{\frac{\beta+ d(\frac1p+\frac1q-1)}{\alpha +\beta +d(\frac1p+\frac1q-1)}}$$ under the conditions on $\beta, p,q$  stated in Theorem    \ref{theo:mikxiao}. We show that the proof of Theorem \ref{theorem:GH} works for any such operator $T$ provided that $\beta<\frac{d}{p'}$ and $q<\infty$.  

First, in light of Remark \ref{remark:general operators}, the upper bounds in Lemma \ref{lemma:local} hold true for $T$, that is, we have
  $$\sup_{f} \frac{\left(\int_{|\xi|\leq t} |Tf|^q\right)^{\frac{1}{q}}}{ \norm{{f}|\xi|^\beta}_p}\lesssim C_{q,p}(t),$$ where $C_{q,p}$ is defined in Lemma \ref{lemma:local}. In particular, for $D_1=D_2=\beta$,
$${\left(\int_{|\xi|\leq t} |Tf|^q\right)^{\frac{1}{q}}} \lesssim t^{\beta + d (\frac{1}{p}+\frac1q -1)} \norm{{f}|\xi|^\beta}_p  .$$
Second, by the triangle inequality and monotonicity,
$$ {\norm{Tf}_q} \lesssim t^{\beta + d (\frac{1}{p}+\frac1q -1)} \norm{{f}|\xi|^\beta}_p  + t^{-\alpha} {\norm{Tf |x|^\alpha}_q}.$$
 Finally, we obtain the result by setting $$t^{\alpha+ \beta +d (\frac1p+\frac1q-1)}=\frac{\norm{Tf |x|^\alpha}_q}{\norm{{f}|\xi|^\beta}_p}.$$
\end{proof}

Second, for the class of special Hadamard operators, it is possible to enlarge the range of parameters for which 
inequality 
\eqref{eq:uncerxiaoprim-} holds true. For instance, we can obtain
\begin{proposition}
       Assume that $q\leq 2$ and either $q\leq p \leq q'$   or $p< \min({q,q'})$.
   Then
$$\sup _T \sup _{ \norm{Tf |x|^\alpha}_q\leq x, \norm{f |\xi|^\beta}_p\leq 1}  {\norm{Tf}_q}\approx H(x),$$ where the supremum is taken over all operators satisfying \eqref{eq:jointboundedness}.

In particular,
 under the conditions on $\beta, p,q$   stated in Theorem
    \ref{theo:mikxiao}, inequality 
    \eqref{eq:uncerxiaoprim-}
    holds for any special Hadamard operator.
\end{proposition}

\begin{proof}
    In light of Proposition \ref{prop:general}, inequality \eqref{eq:uncerxiaoprim-} remains true for all operators satisfying \eqref{eq:jointboundedness} provided that $q<\infty $ and $\beta<\frac{d}{p'}$. If $\beta\geq \frac{d}{p'}$ and $q\leq 2$, letting $d(\frac12-\frac1p)<\gamma<\frac{d}{p'}$ we have
    $$\norm{Tf}^{\alpha+ \theta}_q \lesssim \norm{Tf |x|^{\alpha}}^\theta_q  \norm{f |\xi|^{\gamma}}^{\alpha}_p,$$ where $\theta= \gamma+ d(\frac1p+\frac1q-1)$.
Next, by monotonicity and Lemma \ref{lemma:local} applied to $T^{-1}$ with $D_1=D_2=0$ and our conditions on $p,q$, we have \begin{align*}
    \norm{f |\xi|^{\gamma}}_p &\lesssim t^\gamma \norm{f \mathbbm{1}_{|\xi|\leq t}}_p + t^{\gamma -\beta} \norm{f |\xi|^{\beta}}_p\\
    &\lesssim t^{\gamma+ d(\frac1p+\frac1q-1)} \norm{Tf}_q + t^{\gamma -\beta} \norm{f |\xi|^{\beta}}_p\\
    &\approx \norm{Tf}_q^{\frac{\beta-\gamma}{\theta+ \beta-\gamma}}\norm{f |\xi|^{\beta}}_p^{\frac{\theta}{\theta+ \beta-\gamma}},
\end{align*} where the last equivalence follows by setting $$t^{\beta- \gamma +\theta}= \frac{\norm{f |\xi|^\beta}_p}{\norm{Tf}_q}.$$

In conclusion, we have
     $$\norm{Tf}^{\alpha+ \theta}_q \lesssim \norm{Tf |x|^{\alpha}}^\theta_q  \norm{f |\xi|^{\beta}}^{\frac{\alpha \theta}{\theta+ \beta-\gamma }}_p  \norm{Tf}^{\alpha \frac{\beta-\gamma}{\theta+ \beta -\gamma}}_q,$$ which is equivalent to \eqref{eq:uncerxiaoprim-}.
\end{proof}
Finally, we show that for certain choices of parameters Theorem \ref{theorem:GH}, and in particular inequality 
\eqref{eq:uncerxiaoprim-}, fails for some special Hadamard operators. 
\begin{theorem}
    Let $1\leq p,q\leq \infty$. Then if $\beta,q$ are large enough and $\alpha$ is small enough,
    $$H(x)<\infty$$ but
    $$\sup _T \sup _{ \norm{Tf |x|^\alpha}_q\leq x, \norm{{f} |\xi|^\beta}_p\leq 1}  {\norm{Tf}_q}= \infty,$$ or equivalently,
    $$\sup _T \sup _{ \norm{f |x|^\alpha}_q\leq x, \norm{T{f} |\xi|^\beta}_p\leq 1}  {\norm{f}_q}= \infty,$$where the supremum is taken over all special Hadamard operators.
\end{theorem}
\begin{proof}
For simplicity, we consider $d=1$. Let $N,\varepsilon>0$.
Our first goal is to construct a function $f\in \mathcal{S}$ such that 
\begin{equation}
\label{eq:rearrangement}
    \norm{f}_q \approx N^{1-\frac{1}{q}},  \quad\norm{f|x|^\alpha }_q \approx N^{1-\alpha-\frac{1}{q}},\quad\mbox{ and }\quad  |\widehat{f}|\lesssim \varepsilon^{-1} 
\end{equation}
with the additional property that \begin{equation}
\label{eq:support}
\operatorname{spec}{f} \subset \Sigma:= \bigcup _{j=1}^M [a_j,b_j]  \;\;\mbox{ and } \;\;|\Sigma|\leq 4\varepsilon N. \end{equation}

We use a construction by Nazarov from Section 4.2 in \cite{nazarov1993local}. Let $g:=g_\varepsilon$ be a $1$-periodic $C^\infty$ function satisfying $\int_0^1 g =1$, $g(x)=0$ for $x \in (-1/2,1/2)\setminus (-\varepsilon,\varepsilon)$, and $|g(x)|\lesssim \varepsilon^{-1}$. Let $\psi$ be a smooth function supported on $[-1/2,1/2]$, $\psi_N(\xi)=\psi(\xi/N)$ and set $$F=\widehat{\psi_N}.$$ Define, for $A,b>0$, $$\widehat{f}(\xi)=g_{\varepsilon}(A\xi^2 +b) \psi_N(\xi).$$ Clearly, $|\widehat{f}(\xi)| \lesssim \varepsilon^{-1}$. 

First, we verify \eqref{eq:support}. By definition, $$\operatorname{spec}{f} \subset \Sigma:=  [-N/2,N/2] \cap \bigcup_{k \in \mathbb{Z}} \left\{ \xi: k-\varepsilon \leq A \xi^2+ b \leq k+\varepsilon \right\},$$ from which we see that $\Sigma$ is a finite union of closed intervals. Moreover, for some $0<b<1$ we have that the support of $\Sigma$ has measure $\leq 4\varepsilon N$. Indeed,
$$\int_0^1 |\Sigma(b)| db \leq \sum_{k\in \mathbb{Z}} \int_{\mathbb{R}} \int_{0}^1 \mathbbm{1}_{[-N/2,N/2]}(\xi) \mathbbm{1}_{[k-\varepsilon, k+\varepsilon]} (A\xi^2 +b ) db d \xi \leq 4 \varepsilon N.$$

Second, to show \eqref{eq:rearrangement}, we estimate $|f(x)-F(x)|$ as follows:
\begin{align*}
    |f(x)-F(x)|&=\int_{\mathbb{R}} \widehat{F}(\xi)(1-g_{\varepsilon}(A\xi^2 +b))e^{2 \pi i \xi x} d\xi \\
    &= \int_{\mathbb{R}} \widehat{F}(\xi)\Big(\sum_{k \neq 0} \widehat{g}(k)e^{2 \pi i k(A\xi^2 +b)}º
\Big)e^{2 \pi i \xi x} d\xi\\
     &=\sum_{k \neq 0}  \widehat{g}(k) \int_{\mathbb{R}} \widehat{F}(\xi)
     e^{2 \pi i k(A\xi^2 +b)}e^{2 \pi i \xi x} d\xi.
\end{align*} Observe that, by integration by parts, as obtained in \cite[p. 64]{nazarov1993local},
\begin{align*}
    \left |\int_{\mathbb{R}} \widehat{F}(\xi)e^{2 \pi i k(A\xi^2 +b)}e^{2 \pi i \xi x} d\xi\right| &\lesssim_F \min\Big(\frac{|kA|}{|x|},\frac{1}{|kA|^{\frac12}}\Big),
\end{align*}
and, by the smoothness condition,$$
    |\widehat{g}(k)|\lesssim_g |k|^{-3},\,\, k\in \mathbb{Z}.$$
Thus, \begin{align*}
    |f(x)-F(x)|&\lesssim _{g,F} \sum_{|k| \neq 0} |k|^{-3}\min\Big(\frac{|kA|}{|x|},\frac{1}{|kA|^{\frac12}}\Big)\\
    &\lesssim 
    \begin{cases}
        |A|^{-\frac12}, &|x|\leq |A|^{\frac32},\\
        \frac{|A|}{|x|}&|x|\geq |A|^{\frac32}.
    \end{cases}
    \end{align*}
    Finally,
    $$\norm{ (f -F) |x|^{\alpha}}_q \lesssim_{g,F} |A|^{\frac{3}{2} (\alpha +\frac{1}{q}) -\frac12} + |A|^{\frac{3}{2} (\alpha -1 +\frac{1}{q} )+1}\to 0,\quad |A|\to\infty, $$  provided that $\alpha$ and $\frac{1}{q}$ are small enough.
    Since $\norm{F |x|^{\alpha}}_q\approx  N^{1-\alpha -\frac1q}$ and $\norm{F |x|^{\alpha}}_q\approx  N^{1-\frac1q}$, we obtain \eqref{eq:rearrangement} 
    provided that $|A|$ is sufficiently large.

   After constructing 
$f$, we now proceed to find a suitable special Hadamard operator that will serve as the required counterexample. To this end, we  observe  that if $\sigma:\mathbb{R}\to \mathbb{R}$ is a measure preserving bijection, then $T_\sigma$ defined by $T_\sigma(f)=\widehat{f} (\sigma(\xi))$ is a special Hadamard operator (cf. Example 34 in \cite{dias}). For the function $f$ defined above, let $\sigma$ be a measure preserving bijection such that $\Sigma\subset\sigma([0,5 \varepsilon N])$.
    
    Therefore, if for any function $h$,
    $$\norm{h}_q \lesssim \norm{h |x|^\alpha}_q+ \norm{T_\sigma h |\xi|^\beta}_p $$ holds uniformly on $\sigma$, setting $h=f$ from \eqref{eq:rearrangement} and \eqref{eq:support} we deduce that
  
$$N^{1-\frac{1}{q}} \lesssim  N^{1-\frac{1}{q}-\alpha}+ (\varepsilon N)^{\beta +\frac{1}{p}} \varepsilon^{-1},$$ which is impossible if $\beta>\frac{1}{p'}$ and $\alpha>0$.

    \end{proof}

\section{Appendix. Weighted Fourier inequalities}

In this section, we recall the characterization of the radial non-increasing weights $u,v^{-1}$  for which the weighted Fourier inequality 
\begin{align}
\label{ineq:pittgen}
\big\|u\widehat{f}\,\big\|
_{q}&\leq C \norm{ f v}_{p}
\end{align}
holds. For further discussion of the historical development of this topic, see \cite{benedetto2003weighted} and \cite{saucedo2025fouriertransformextremizerclass}.

\begin{theorem}[Theorem 5.1 in \cite{saucedo2025fouriertransformextremizerclass}]
\label{theorem:mainh}
    Let $1\leq p \leq \infty$ and $0<q\leq \infty$. Set $\frac{1}{r}=\frac{1}{q}-\frac{1}{p},$ $ \frac{1}{p^\sharp}=\frac{1}{2}-\frac{1}{p}$, and $\frac{1}{q^\sharp}= \frac{1}{q}-\frac{1}{2}$. Assume that $u$ and $v^{-1}$ are radial non-increasing weights. Define     \begin{equation*}
     U(t)=\int_0 ^t u^{*,q},\quad V(t)= \int_0 ^t v_*^{p},\quad\mbox{and}\quad
     \xi(t) = U(t)+ t^{\frac{q}{2}} \left(\int_t ^\infty u^{*,q^{\sharp}} \right)^{\frac{q}{q^{\sharp}}}.\end{equation*}
     Then
     inequality \eqref{ineq:pittgen} holds if and only if
    \begin{enumerate}[label=(\roman*)]
    \item $p \leq q$ and  \begin{equation}\label{condition:11}
        C_1:=\sup _{s>0} \left(\int_0^s u^{*,q} \right)^{\frac{1}{q}} \left(\int_0^{1/s} v_*^{-p' } \right)^{\frac{1}{p'}} <\infty; 
        \end{equation}
        \item $1\leq q<p$, $2\leq\max(q,p'),$ and  \begin{equation}\label{nec11}
        \begin{aligned}
      C_2:&= \Bigg(
       \int_0^\infty u^{*,q}(s)
       \Big(\int_{0}^{s} u^{*,q}\Big)^{\frac{r}{p}}\Big(\int_{0}^{1/s} v_*^{-p'}\Big)^{\frac{r}{p'}}  ds \Bigg)^{\frac{1}{r}}\\
       &\approx \Bigg(
       \int_0^\infty v_*^{-p'}(s)
       \Big(\int_{0}^{s} v_*^{-p'}\Big)^{\frac{r}{q'}}\Big(\int_{0}^{1/s} u^{*,q}\Big)^{\frac{r}{q}}  ds \Bigg)^{\frac{1}{r}}<\infty;
       \end{aligned}
    \end{equation}

     \item $q<2< p<\infty$, 
     $\int_1 ^\infty u^{*,q^{\sharp}} < \infty
     $
     and  
$C_3=C_2+C_4<\infty,
$ 
 where 
     \begin{eqnarray*}
C_4&=&\Bigg(\int_0^\infty u^{*,q}(t) U^{\frac{q^{\sharp}}{2}}(t) \xi^{-\frac{q^{\sharp}}{2}}(t) t^{-\frac{q}{2}} \Big(\int_{t}^ \infty u^{*,q}(r) U^{\frac{q^{\sharp}}{2}}(r) \xi^{-\frac{q^{\sharp}}{2}}(r) r^{-\frac{q}{2}}   dr \Big)^{\frac{r}{p}}
\\&&\qquad \times \Big( \int_{t^{-1}} ^\infty v_*^{-p^{\sharp}}(r) \Big(\frac{V(r)}{r v_*^{p}(r)}\Big)^{-\frac{p^{\sharp}}{2}} dr \Big) ^{\frac{r}{p^{\sharp}}}dt\Bigg)^{\frac{1}{r}};
     \end{eqnarray*}
 \item $q<2< p=\infty$, 
     $\displaystyle
         \int_1 ^\infty u^{*,q^{\sharp}} < \infty
    $
      and  
     \begin{equation}
   \label{eq:pinfty}
   C_5:= \left(\int_0 ^\infty U^{\frac{q^{\sharp}}{2}}(t) u^{*,q}(t) \xi^{-\frac{q^{\sharp}}{2}}(t) t^{-\frac{q}{2}} \Big(\int_0 ^t \Big ( \int_0 ^{r^{-1}} v_*^{-1} \Big)^2 dr \Big)^{\frac{q}{2}}dt\right)^{\frac{1}{q}}<\infty;
\end{equation}

\item $q<1\leq p\leq 2$,
$\displaystyle
         \int_1 ^\infty u^{*,q^{\sharp}} < \infty,
     $
     and  
     $
C_6=C_2+C_{7}<\infty,
$

\end{enumerate}
    where
     \begin{eqnarray*}
C_{7}&=&\Bigg(\int_0^\infty u^{*,q}(t) U^{\frac{q^{\sharp}}{2}}(t) \xi^{-\frac{q^{\sharp}}{2}}(t) t^{-\frac{q}{2}} 
\\&&\qquad \times\Big(\int_{t}^ \infty u^{*,q}(r) U^{\frac{q^{\sharp}}{2}}(r) \xi^{-\frac{q^{\sharp}}{2}}(r) r^{-\frac{q}{2}}   dr \Big)^{\frac{r}{p}} \Big(\sup_{1/t\leq y<\infty} \frac{y ^{\frac{r}{2}}}{V^{\frac{r}{p}}(y)}
\Big)dt\Bigg)^{\frac{1}{r}}.
     \end{eqnarray*}     
Moreover, the optimal constant $C$ in the 
inequality $
 \big\|u\widehat{f}\,\big\|_{q}\leq C \norm{ f v}_p
 $
 satisfies
$$C\approx C_j, \qquad j=1,2,3,5,6,
$$
in items (i)-(v), respectively.
\end{theorem}
\begin{remark}
\label{remark:general operators}
    For general linear operators $T$ satisfying $\norm{Tf}_2 \lesssim \norm{f}_2$ and $\norm{Tf}_\infty\lesssim \norm{f}_1$,  the optimal constant $C$ in the 
inequality $
 \big\|uTf\,\big\|_{q}\leq C \norm{ f v}_p
 $
 satisfies
$$C\lesssim C_j, \qquad j=1,2,3,5,6,
$$
in items (i)-(v), respectively.
\end{remark}
In this work  we  only use items $(i)$ and $(ii)$. To deal with the case $q<2<p$, instead of the more difficult case $(iii)$ it suffices to use a simpler necessary condition.
\begin{lemma}[Lemma 3.2 in \cite{saucedo2025fouriertransformextremizerclass}]
\label{lemma:prelinec}
    Let $2< p\leq \infty$ and $0<q\leq \infty$. Assume that inequality \eqref{ineq:pittgen} holds for $f\in L^1$. Then
\begin{equation}\label{new-n}
 \sup_{s} \left( \int_{|x|\leq s^{-1}} u^q  \right)^{\frac{1}{q}}   \Bigg(\sum_{n \in \mathbb{Z}^d} \Big(\int_{sn+ s[-\frac{1}{2}, \frac{1}{2}]^d}v^{-p'}\Big) ^{\frac{p^{\sharp}}{p'}} \Bigg) ^{\frac{1}{p^{\sharp}}} \lesssim C.
 \end{equation}
In particular, if $v$ is radial, non-decreasing,

\begin{equation}\label{new-n2}
 \sup_{s}  s^{\frac{d}{2}}\Big( \int_{|x|\leq s^{-1}} u^q  \Big)^{\frac{1}{q}}   \Big(\int_{|x|\geq s} v^{-p^\sharp}\Big)^{\frac{1}{p^\sharp}} \lesssim C.
 \end{equation}
    
\end{lemma}

\bibliographystyle{abbrv}
\bibliography{main.bib}
\end{document}